\documentclass[12pt,english,a4paper]{smfart}

\usepackage[french]{babel}
\usepackage{amssymb}
\usepackage{latexsym}
\usepackage{bm}

\numberwithin{equation}{part}

\def\di{\displaystyle}

\def\E{\cal E}
\def\Ef{{}_a^{\alpha} {\bf E}_{b}^{\beta}}
\def\p{{\bf p}}
\def\L{\mbox{\rm L}}
\def\Ll{\cal L}
\def\Llf{{\cal L}_{a,b}^{\alpha ,\beta}}

\def\M{\mbox{\rm M}}

\def\N{\mbox{\rm N}}
\def\D{{\cal D}^{\alpha ,\beta}_{\mu}}
\def\Di{{\cal D}^{\beta ,\alpha}_{-\mu}}
\def\RL{\mbox{\rm RL}}
\def\adta{{}_a \mbox{\rm D}_t^{\alpha}}
\def\tdba{{}_t \mbox{\rm D}_b^{\alpha}}
\def\adtb{{}_a \mbox{\rm D}_t^{\beta}}
\def\tdbb{{}_t \mbox{\rm D}_b^{\beta}}
\def\adtm{{}_a \mbox{\rm D}_t^m}
\def\tdbm{{}_t \mbox{\rm D}_b^m}
\def\aita{{}_a \mbox{\rm D}_t^{-\alpha}}
\def\tiba{{}_t \mbox{\rm D}_b^{-\alpha}}
\def\lap{\mbox{\rm LAP}}
\def\flap{\mbox{\rm FLAP}}
\def\Rev{\mbox{\rm Rev}^{\alpha}}
\def\Revo{\mbox{\rm\bf Rev}}
\def\FD{\mathbb{D}}
\def\EL{\mbox{\rm EL}}
\def\FEL{\mbox{\rm FEL}}
\def\Op{\mbox{\rm O}}

\def\F{\mbox{\rm F}}

\def\R{\mathbb{R}}
\def\N{\mathbb{N}}

\def\C{\mathbb{C}}

\def\P{\mathbb{P}}
\def\ds{\mbox{\em ds}}
\def\Emb{\di {}_{a}^{\alpha} \mbox{\rm Emb}_{b}^{\beta} (\mu)}
\def\Embi{\di {}_{a}^{\beta} \mbox{\rm Emb}_{b}^{\alpha} (-\mu)}
\def\var{\mbox{\rm\bf Var}(a,b)}
\def\Var{\mbox{\rm\bf Var}}
\def\vari{\mbox{\rm\bf Var}(a,b)}
\def\varf{{}_a^{\alpha} \mbox{\rm\bf Var}^{\beta}_b}
\def\Id{\mbox{\rm Id}}

\def\cal{\mathcal}

\usepackage{amsmath}
\usepackage{amssymb}
\usepackage{amsthm}
\usepackage{a4wide}

\newtheorem{thm}{Theorem}[part]
\newtheorem{defi}{Definition}[part]
\newtheorem{lem}{Lemma}[part]

\newtheorem{cor}{Corollary}[part]

\newtheorem{rema}{Remark}[part]

\title[Fractional embedding]{Fractional embedding of differential operators and Lagrangian systems}
\author{Jacky CRESSON}
\email{jacky.cresson@univ-pau.fr} \email{cresson@ihes.fr}

\address{Universit\'e de Pau et des Pays de l'Adour, Laboratoire de Math\'ematiques appliqu\'ees de Pau, CNRS UMR 5142}
\address{I.H.E.S, 35 route de Chartres, 91440 Bures sur Yvette.}

\begin{document}
\baselineskip 6mm
\maketitle

\setcounter{tocdepth}{3}

\begin{abstract}
This paper is a contribution to the general program of embedding theories of dynamical systems. Following our previous work on the Stochastic embedding
theory developed with S. Darses, we define the fractional embedding of differential operators and ordinary differential equations. We construct an
operator combining in a symmetric way the left and right (Riemann-Liouville) fractional derivatives. For Lagrangian systems, our method provide a
fractional Euler-Lagrange equation. We prove, developing the corresponding fractional calculus of variations, that such equation can be derived via
a fractional least-action principle. We then obtain naturally a fractional Noether theorem and a fractional Hamiltonian formulation of fractional
Lagrangian systems. All these constructions are coherents, {\it i.e.} that the embedding procedure is compatible with the fractional calculus of variations. We
then extend our results to cover the Ostrogradski formalism. Using the fractional embedding and following a previous work of F. Riewe, we obtain a fractional Ostrogradski formalism which allows
us to derive non-conservative dynamical systems via a fractional generalized least-action principle. We also discuss the Whittaker equation and obtain a
fractional Lagrangian formulation. Last, we discuss the fractional embedding of continuous Lagrangian systems. In particular, we obtain a fractional
Lagrangian formulation of the classical fractional wave equation introduced by Schneider and Wyss as well as the fractional diffusion equation.
\end{abstract}


\begin{tiny}
\tableofcontents
\end{tiny}

\newpage
\begin{verse}
\tiny
\hfill "Il y a de l'apparence qu'on tirera un jour\\
\hfill des consequences bien utiles de ces paradoxes\footnote{Derivatives of non-integer order},\\
\hfill car il n'y a gueres de paradoxes sans utilit\'e."\footnote{It will lead to a paradox, from which one day useful consequences will be drawn.}
\vskip 2mm
\hfill Leibniz, Letter to L`Hospital, September 30, 1695
\end{verse}
\vskip 5mm
\begin{verse}
\tiny
\hfill Comme la construction du monde est la plus parfaite possible\\
\hfill et qu'elle est due \`a un Cr\'eateur infiniment sage,\\
\hfill il n'arrive rien dans le monde qui ne pr\'esente quelques propri\'et\'e de maximum ou minimum. \\
\hfill Aussi ne peut-on douter qu'il soit possible de d\'eterminer tous les effets de
l'Univers\\
\hfill par leurs causes finales, \`a l'aide de la m\'ethode des maxima et minima\\
\hfill avec tout autant de succ\`es que par leurs causes efficientes.
\vskip 3mm
\hfill {\it Methodus inveniendi  lineas curveas maximi  minimive proprietate gaudentes}

\hfill in Eulerii Opera omnia, S\'erie I, 24,

\hfill Berlin-Basel-Boston-Stuttgart, Lipsiae-Birkhauser Verlag, 1911
\end{verse}

\part*{Introduction}

The aim of this paper is to introduce a general procedure called the {\it fractional embedding procedure}, which roughly speaking allows us
to associate a fractional analogue of a given ordinary differential equation in a more or less {\it canonical} way. The fractional embedding
procedure is part of a global point of view on dynamical systems called {\it embedding theories of dynamical systems} \cite{cr2}. Before
describing the general strategy underlying all embedding theories and the particular fractional embedding procedure, we provide a set of
problems which lead us to our point of view:\\

- {\it Turbulence}: Fluid dynamics is modelled by partial differential equations. Solutions of these equations must be sufficiently smooth. However,
there exits {\it turbulent} behavior which correspond to very irregular trajectories. If the underlying equation has a physical meaning, then one must give
a sense to this equation on irregular functions. This remark is the starting point of Jean Leray's work on fluid mechanics \cite{ler}. He introduces what he
calls {\it quasi-derivation} and the notion of weak-solutions for PDE. This first work has a long history and descendence going trough the definition of Laurent
Schwartz's {\it distribution} and {\it Sobolev spaces}. We refer to \cite{and} for an overview of Jean Leray's work in this domain.\\

- {\it Deformation Quantization problems}: The problem here is to go from classical mechanics to quantum mechanics trough a deformation involving
the Planck constant. Roughly speaking, we have a one parameter $h$ family of
spaces and operators such that they reduce to usual spaces and operators when $h$ goes to zero. For example, we can look for a deformation of
the classical derivative using its algebraic characterization trough the Leibniz rule. Another way is following L. Nottale \cite{no} to assume that the space-time at the atomic
scale is a non-differentiable manifold. In that case, we obtain a one parameter smooth deformation of space-time by smoothing at different scales.
The main problem is then to look for the deformation of the classical derivative during this process. We refer to \cite{cr3} for more details.\\

- {\it Long term behavior of the Solar-system}: The dynamics of the Solar system is usually modelled by a $n$-body problem. However, the study of
the long-term behavior must include several perturbations terms, like tidal effects, perturbations due to the oblatness of the sun, general relativity
effects etc. We do not know the whole set of perturbations which can be of importance for the long term dynamics. In particular, it is not clear that
the remaining perturbations can be modelled using ordinary differential equations. Most of stability results uses in the Solar systems dynamics make
this assumption implicitly \cite{marmi}. An idea is to try to look for the dynamics of the initial equation on more general objects like stochastic processes, by
extending the ordinary derivative. Then one can look for the stochastic perturbation of the underlying stochastic equation which contains the original one.
As a consequence, we can provide a set of dynamical behaviors which have a strong significance being stable under very general perturbation terms. This
strategy is developed in (\cite{cd1} \cite{cd3}) and applied in \cite{cd4}.\\

These problems although completely distinct have a {\it common core}: we need to {\it extend} the classical derivative to a more general functional space. This
extension being given, we have a natural, but not always {\it canonical}, associated equation.\\

However, we are lead to two completely distinct theories depending on the nature of the {\it extension} we make for the classical derivative.\\

For the Solar-system
problem, we need to extend the classical derivative to stochastic processes by imposing that the new operator reduces to the classical derivative on
differentiable deterministic processes. In that case, the initial equation is present in the extended one and we use the terminology of {\it embedding theory}.
We describe the strategy with more details in the next section.\\

For quantization problems, in particular first quantization of classical mechanics, we
have an extra parameter $h$ (the Planck constant) and the extended operator denoted by ${\cal D}_h$ reduces to the classical derivative when $h=0$. The
initial equation is not contain in the extended one, but we have a {\it continuous} deformation of this equation depending on $h$. We then use the terminology of
{\it deformation theory}. We give more details on this type of theories in the following.\\

\section{Embedding theories}

The general scheme underlying embedding theories of ordinary or partial differential equations is the following:
\vskip 3mm
\begin{itemize}
\item Fix a functional space $\cal F$ and a mapping $\iota :C^0 \rightarrow {\cal F}$.
\item Extend the ordinary derivative on $\cal F$ by imposing a gluing rule.
\item Extend differential operators.
\item Extend ordinary or partial differential equations.
\end{itemize}
\vskip 3mm
Let us denote by $\cal D$ the extended derivative on $\cal F$. The gluing rule impose that the extended derivative reduces to the ordinary derivative on
$\iota (C^1 )$, {\it i.e.} that we have  ${\cal D} \iota (x)=\iota (\dot{x} )$, for all $x\in C^1$, where $\dot{x}=dx/dt$. As a consequence, the
original equation can be recovered via the embedded equation by restricting the underlying functional space to $\iota (C^k )$, $k$ depending on the
order of the original equation.\\

Most of the time, we have in mind applications to physics where Lagrangian systems play a fundamental role. The importance of these systems is related to
the fact that they can be derived via a {\it first principle}, the least-action principle. Moreover, in some cases of importance we can find
a symmetric representation of these systems as {\it Hamiltonian systems}. Hamiltonian systems are fundamental by many aspects. Of particular interest is
the fact that they can be {\it quantized} in order to obtain {\it quantum} analogues of classical dynamical systems.\\

An embedding procedure can always be applied to Lagrangian systems. We obtain to embed objects which are the {\it embed action functional} and the
{\it embed Euler-Lagrange equation}. At this point, an embedding procedure can be considered as a particular quantization procedure related to
the underlying functional space. However, we have a more elaborate picture here. As we have an embed action functional, we can develop the
associated calculus of variations, that we call the embed calculus of variations in what follows. We then obtain an embed least-action
principle with an associated Euler-Lagrange equation. As a consequence, we have two kinds of embed Euler-Lagrange equations and it is not clear that the
one obtain directly from the embedding of the classical Euler-Lagrange equation and the one obtain via the embed least-action principle coincide. This
problem is recurrent in all embedding theories of Lagrangian systems and is called the {\it coherence problem}. This problem is far from
being trivial in most of the already existing embedding theories, like the stochastic one (\cite{cd1} \cite{cd3}) or the quantum one \cite{cr4}.

\section{Emergence of fractional derivatives}

In this paper, we develop an embedding theory of ordinary differential equations and Lagrangian systems using {\it fractional derivatives}. Precise
definitions will be given in section \ref{part1}. We only give an heuristic introduction for these operators and some basic problems where they arise
naturally.\\

A fractional derivative is an operator which gives a sense to a real power of the classical differential operator $d/dt$, {\it i.e.} that we want to
consider an expression like
\begin{equation}
\di {d^{\alpha}\over dt^{\alpha}} , \ \ \alpha \geq 0 .
\end{equation}
The previous problem appears for the first time in a letter from Leibniz to L'H\^opital in 1695: "{\it Can the meaning of derivatives with integer order
be generalized to derivatives with non-integer orders} ?". Many mathematicians have contributed to this topic including Leibniz, Liouville, Riemann ,
etc. We refer to \cite{os} or \cite{mr} for a historical survey. A number of definitions have emerged over the years including
Riemann-Liouville fractional derivative, Grunwald-Letnikov fractional derivative, Caputo fractional derivative, etc. In this article we restrict our
attention to the Riemann-Liouville fractional derivative, although the embedding theory can be developed for an arbitrary given fractional calculus with
different technical difficulties.\\

The main difficulties when dealing with fractional derivatives are related to the following properties:\\

(i) fractional differential operators are not {\it local} operators

(ii) the adjoint of a fractional differential operator is not the negative of itself\\

Property (i) is widely use in applications and explain part of the interest for these operators to model phenomenon with {\it long memory} (see for
example \cite{comte}).\\

Other problems arise during computations. Developing the fractional calculus of variation and the associated results (the fractional Noether theorem)
we have encountered difficulties linked with the following facts:\\

(i) the classical {\it Leibniz rule} $(fg)'=f'g+fg'$ is more complicated (see \cite{skm})

(ii) there exists no simple formula for the fractional analogue of the {\it chain rule} in classical differential calculus (see \cite{skm})\\

This last difficulty is of special importance in the derivation of the fractional Noether theorem.\\

The fractional framework has been used in a wide variety of problems. We note in particular applications in {\it turbulence} \cite{clz},
{\it chaotic dynamics} \cite{zu} and {\it quantization} \cite{muba}.\\

In this paper, we will frequently quote the work of F. Riewe (\cite{ri1},\cite{ri2}) which proposes a fractional approach to {\it nonconservative} dynamical
systems. The main property of these systems is that they induce an {\it arrow of time} due to irreversible dissipative effects. The relation between
fractional derivatives, nonconservative systems and irreversibility have been discussed for example in the book \cite{mnn}.\\

Irreversibility implies that we look for the past ${\cal P}_t$ and the future ${\cal F}_t$ of a given given dynamical process $x(s)$, $s\in \R$ at time $t$,
{\it i.e.} on the information ${\cal P}_t =\{ x(s), a\leq s \leq t\}$ and ${\cal F}_t =\{ x(s) ,\ t\leq s\leq b \}$ where $a$ and $b$ can be chosen and
depends on the amount of information we are keeping from the past and the future. This induce the fact that we look for two quantities, not yet defined
that we denote by $d_- x(t)$ and $d_+ x(t)$ from the point of view of derivatives.\\

The past and future information can be {\it weighted}, {\it i.e.} that
we look not for $x(s)$ but to $w(s,t) x(s)$ where $w(s,t)$ give the importance of the information at time $s$ with respect to time $t$. This can be
achieved using a weight $1\over \mid t-s \mid^{\alpha +1}$ and regularizing the corresponding function. We then are lead to two quantities $d_-^{\alpha} x(t)$
and $d_+^{\alpha} x(t)$, which represent a weighted information on the past and future behavior of the dynamical process.\\

The previous idea is well formalized by the left and right (Riemann-Liouville) derivatives \cite{skm}. We refer to part \ref{part1} for precise
definitions.

\section{Deformation theories and the fractional framework}

The fractional framework follows the general strategy outline in the previous section. However, a new ingredient comes into play which makes the
fractional embedding different from the existing stochastic or quantum embedding theories. We have used in this paper the {\it left and right Riemann-Liouville
derivatives} with different indices for the left and right differentiation, {\it i.e.} we consider $\adta$ and $\tdbb$. The extended operator depends
naturally on these two operators and is denoted ${\cal D}^{\alpha ,\beta}$. However, this operator does not reduce to the ordinary derivative on the set
of differentiable functions. We recover the ordinary derivative only when $\alpha=\beta=1$. As a consequence, we can associate to a given ordinary
differential equation a two parameters family of {\it fractional differential equations}. The original equation being recovered for a special choice of
these parameters, {\it i.e.} $\alpha=\beta=1$. In that case, we propose to the use the terminology of {\it fractional deformation} and to keep
the terminology of {\it fractional embedding} for the procedure which associated a fractional analogue of an ordinary differential equation.\\

Deformation theory can be formalized as follows:\\

A {\it deformation theory} is the data of:
\vskip 3mm
\begin{itemize}
\item A {\it finite} set ${\cal P}=\{ (p_1 ,\dots ,p_{\nu}),\ p_i \in A \}$ of parameters where $A$ is a given interval of $\R$.
\item A $\nu$-parameter family of functional spaces ${\cal F} =\{ {\cal F}_{\cal P} \}_{{\cal P}\in A^{\nu} }$.
\item Operators ${\cal D}_{\cal P}$ defined on ${\cal F}_{\cal P}$ such that there exists ${\cal P}_0 \in A^{\nu}$ satisfying $C^1 \subset
{\cal F}_{{\cal P}_0}$ and ${\cal D}_{{\cal P}_0} (x)=\dot{x}$ for $x\in C^1$.
\end{itemize}

The condition on $A$ is only here to be sure that we have a {\it continuous} dependance of the whole construction on the parameters.\\

The main difference between deformation and embedding lies in the fact that it is no usually easy to obtain {\it information} on the initial equation
from the deformed one. We must use {\it asymptotic methods}, looking for the behavior of the deformed equation when $p\rightarrow p_0$. This is not the
case for a true embedding theory as the initial equation is already present in the embedded one.

\section{Plan of the paper}

Our paper has the same architecture as our previous monograph \cite{cd3} with S\'ebastien Darses about the stochastic embedding of dynamical systems. As
a consequence, the comparison between the two embedding procedure will be easier.\\

In part \ref{part1} we recall the definitions of the left and right fractional derivatives. We also define left and right fractional derivatives which
have satisfy a semi-group property and the adequate functional spaces on which they are defined following a previous work of
Erwin and Roop \cite{er}. We also recall a product rule formula.\\

In part \ref{part2} we define the fractional embedding of differential operators and ordinary differential equations.\\

In part \ref{part3} we study the fractional embedding of Lagrangian systems. We obtain a fractional analogue of the Euler-Lagrange equations.\\

In part \ref{fraccalculvar} we develop a fractional calculus of variations associated to the fractional embedding of classical functionals. generalizing
a previous work of O.P Agrawal \cite{agra} We prove
two versions of the least action principle depending on the underlying authorized space of variations. We prove in particular a coherence theorem,
which roughly speaking state that the fractional embedding of the Euler-Lagrange equation coincide with the fractional Euler-Lagrange equation
obtained via the fractional least-action principle.\\

Part \ref{part5} study the behavior of symmetries under the fractional embedding procedure. In particular, we prove a fractional Noether theorem which
generalizes a recent result of Frederico and Torres \cite{gasto}.\\

In part \ref{part6} we derive the analogue of the Hamilton formalism for our fractional Lagrangian systems.\\

In part \ref{part7} we extend results of parts \ref{part3} and \ref{fraccalculvar} to cover the Ostrogradski formalism for Lagrangian systems. In this unified
framework, we recover classical results of F. Riewe (\cite{ri1} \cite{ri2}). Precisely, we obtain a fractional Lagrangian derivation of Nonconservative
systems.\\

In part \ref{part8} we study the fractional embedding of continuous Lagrangian systems. In particular, we prove that the classical fractional wave
equation introduced by W.R. Schneider and W. Wyss \cite{sw} under ad-hoc assumptions, is the fractional embedding of the classical wave equation which
respects the underlying continuous Lagrangian structure of the equation. An analogous result is obtained for the fractional diffusion equation.\\

We then conclude with some open problems and perspectives.

\newpage
\part{Fractional operators}
\label{part1}
\setcounter{section}{0}

In \cite{agra} Agrawal has studied Fractional variational problems using the Riemann-Liouville derivatives. He notes that even if the initial functional
problems only deals with the left Riemann-Liouville derivative, the right Riemann-Liouville derivative appears naturally during the computations. In this
section, we construct an operator combining the left and right Riemann-Liouville (RL) derivative. We remind some results concerning functional spaces
associated to the left and right RL derivative. In particular, we discuss the possibility to obtain a law of exponents.

\section{Fractional differential operators}

\subsection{Left and Right Riemann-Liouville derivatives}

We define the left and right Riemann-Liouville derivatives following (\cite{os} \cite{skm} \cite{pod} \cite{mr}).

\begin{defi}[Left Riemann-Liouville Fractional integral] Let $x$ be a function defined on $(a,b)$, and $\alpha >0$. Then the left Riemann-Liouville
fractional integral of order $\alpha$ is defined to be
\begin{equation}
\aita x(t):=\di {1\over \Gamma (\alpha )} \di\int_a^t (t-s)^{\alpha -1} x(s) ds .
\end{equation}
\end{defi}

\begin{defi}[Right Riemann-Liouville Fractional integral] Let $x$ be a function defined on $(a,b)$, and $\alpha >0$. Then the right Riemann-Liouville
fractional integral of order $\alpha$ is defined to be
\begin{equation}
\tiba x(t):=\di {1\over \Gamma (\alpha )} \di\int_t^b (s-t)^{\alpha -1} x(s) ds .
\end{equation}
\end{defi}

Left and right (RL) integrals satisfy some important properties like the semi-group property. We refer to \cite{skm} for more details.

\begin{defi}[Left and Right Riemann-Liouville fractional derivative]
Let $\alpha >0$, the left and right Riemann-Liouville derivative of order $\alpha$, denoted by $\adta$ and $\tdba$ respectively, are defined by
\begin{equation}
\adta x(t) =\di {1\over \Gamma (n-\alpha )} \di \left ( \di {d\over dt} \right )^n \int_a^t (t-s)^{n-\alpha-1} x(s) ds ,
\end{equation}
and
\begin{equation}
\tdba x(t) =\di {1\over \Gamma (n-\alpha )} \di \left ( -\di {d\over dt} \right )^n \int_t^b (t-s)^{n-\alpha-1} x(s) ds ,
\end{equation}
where $n$ is such that $n-1 \leq \alpha <n$.
\end{defi}

If $\alpha=m$, $m\in \N^*$, and $x\in C^m (]a,b[)$ we have
\begin{equation}
\adtm x=\di {d^m x\over dt^m}  ,\ \ \tdbm =\di -{d^m x\over dt^m} .
\end{equation}
This last relation which ensures the gluing of the left and right Riemann-Liouville ($\RL$) derivative to the classical derivative will be of
fundamental importance in what follows.\\

If $x(t) \in C^0$ with left and right-derivatives at point $t$ denoted by $\di {d^+ x\over dt}$ and $\di {d^- x\over dt}$ respectively then
\begin{equation}
\adtm x =\di {d^+ x\over dt}  ,\ \ \tdbm =\di {d^- x\over dt} .
\end{equation}

In what follows, we denote by ${}_a^{\alpha} \mbox{\rm\bf E}$, $\mbox{\rm\bf E}_b^{\beta}$ and $\Ef$ the functional spaces defined by
\begin{equation}
{}_a^{\alpha} \mbox{\rm\bf E} =\{ x\in C([a,b]),\ \adta x\ \mbox{\rm exists} \} ,\ \ \
\mbox{\rm\bf E}_b^{\beta} =\{ x\in C([a,b]),\ \tdbb x\ \mbox{\rm exists} \} ,
\end{equation}
and
\begin{equation}
\Ef ={}_a^{\alpha} \mbox{\rm\bf E} \cap \mbox{\rm\bf E}_b^{\beta} .
\end{equation}

\begin{rema}
Of course the set $\Ef$ is non-empty. Following (\cite{skm} Lemma 2.2 p.35) we have $AC([a,b]) \subset \Ef$, where $AC([a,b])$ is the set of
absolutely continuous functions on the interval $[a,b]$ (see \cite{skm} Definition 1.2).
\end{rema}

The operators of ordinary differentiation of integer order satisfy a commutativity property and the law of exponents (the semi-group
property) {\it i.e.}
\begin{equation}
\di {d^n \over dt^n} \circ \di {d^m \over dt^m} = \di {d^m \over dt^m} \circ \di {d^n \over dt^n}=\di {d^{n+m} \over dt^{n+m}} .
\end{equation}
These two properties in general fail to be satisfied  by the left and right fractional RL derivatives. We refer to (\cite{mr} $\S$.IV.6) and (\cite{gm} p.233) for more details
and examples. These bad properties are responsible for several difficulties in the study of fractional differential equations. We refer to \cite{pod}
for more details.

\subsection{Left and right fractional derivatives}

In some cases, we need that our fractional operators satisfy additional properties like the {\it semi-group property}. Following \cite{er} we introduce the {\it left
and right fractional derivatives} as well as convenient functional spaces on which we have the semi-group property.

\begin{defi}[Left fractional derivative] Let $x$ be a function defined on $\R$, $\alpha >0$, $n$ be the smallest integer greater than $\alpha$ ($n-1
\leq \alpha <n$), and $\sigma =n-\alpha$. Then the left fractional derivative of order $\alpha$ is defined to be
\begin{equation}
\FD^{\alpha} x (t):={}_{\infty} D_t^{\alpha} x (t)=\di {d^n \over dt^n} {}_{\infty} D_t^{-\alpha} x(t)=\di {1\over \Gamma (\sigma) }\di {d^n\over dt^n}
\int_{-\infty}^t (t-s)^{\sigma -1} x(s) ds .
\end{equation}
\end{defi}

\begin{defi}[Right fractional derivative] Let $x$ be a function defined on $\R$, $\alpha >0$, $n$ be the smallest integer greater than $\alpha$ ($n-1
\leq \alpha <n$), and $\sigma =n-\alpha$. Then the right fractional derivative of order $\alpha$ is defined to be
\begin{equation}
\FD_*^{\alpha} x (t):={}_t D_t^{\infty} x (t)=\di (-1)^n {d^n \over dt^n} {}_t D_{\infty}^{-\alpha} x(t)=\di {(-1)^n\over \Gamma (\sigma) }\di {d^n\over dt^n}
\int_t^{\infty} (s-t)^{\sigma -1} x(s) ds .
\end{equation}
\end{defi}

If $\overline{\mbox{\rm Supp} (x)} \subset (a,b)$ we have $\FD^{\alpha} x =\adta x$ and $\FD_*^{\alpha} x =\tdba x$.\\

In \cite{er} several useful functional spaces are introduced. Let $I\subset \R$ be an open interval (which may be unbounded). We denote by
$C_0^{\infty} (I)$ the set of all functions $x\in C^{\infty} (I)$ that vanish outside a compact subset $K$ of $I$.

\begin{defi}[Left fractional derivative space] Let $\alpha >0$. Define the semi-norm
\begin{equation}
\mid x \mid_{J_L^{\alpha} (\R )} :=\parallel \FD^{\alpha} x \parallel_{L^2 (\R)} ,
\end{equation}
and norm
\begin{equation}
\parallel x\parallel_{J_L^{\alpha} (\R )} :=\left ( \parallel x\parallel^2_{L^2 (\R )} +\mid x \mid^2_{J_L^{\alpha} (\R )} \right ) ^{1/2} .
\end{equation}
and let $J_L^{\alpha} (\R )$ denote the closure of $C_0^{\infty} (\R )$ with respect to $\parallel \cdot \parallel_{J_L^{\alpha} (\R )}$.
\end{defi}

Similarly, we can defined the right fractional derivative space.

\begin{defi}[Right fractional derivative space] Let $\alpha >0$. Define the semi-norm
\begin{equation}
\mid x \mid_{J_R^{\alpha} (\R )} :=\parallel \FD_*^{\alpha} x \parallel_{L^2 (\R)} ,
\end{equation}
and norm
\begin{equation}
\parallel x\parallel_{J_R^{\alpha} (\R )} :=\left ( \parallel x\parallel^2_{L^2 (\R )} +\mid x \mid^2_{J_R^{\alpha} (\R )} \right ) ^{1/2} .
\end{equation}
and let $J_R^{\alpha} (\R )$ denote the closure of $C_0^{\infty} (\R )$ with respect to $\parallel \cdot \parallel_{J_R^{\alpha} (\R )}$.
\end{defi}

We now assume that $I$ is a {\it bounded} open subinterval of $\R$. We restrict the fractional derivative spaces to $I$.

\begin{defi}
Define the spaces $J_{L,0}^{\alpha} (I)$, $J_{R,0}^{\alpha} (I)$ as the closure of $C_0^{\infty} (I)$ under their respective norms.
\end{defi}

These spaces have very interesting properties with respect to $\FD$ and $\FD_*$. In particular, we have the following semi-group property:

\begin{lem}
\label{semigroup}
For $x\in J_{L,0}^{\beta} (I)$, $0<\alpha <\beta$ we have
\begin{equation}
\FD^{\beta} x =\FD^{\alpha} \FD^{\beta -\alpha} x
\end{equation}
and similarly for $x\in J_{R,0}^{\beta} (I)$,
\begin{equation}
\FD_*^{\beta} x =\FD_*^{\alpha} \FD_*^{\beta -\alpha} x .
\end{equation}
\end{lem}

We refer to (\cite{er} Lemma 2.9) for a proof.\\

The fractional derivative spaces $J_{L,0}^{\alpha} (I)$ and $J_{R,0}^{\alpha} (I)$ have been characterized when $\alpha >0$. We denote by
$H_0^{\alpha} (I)$ the {\it fractional Sobolev space}.

\begin{thm}
Let $\alpha >0$. Then the $J_{L,0}^{\alpha} (I)$, $J_{R,0}^{\alpha} (I)$ and $H_0^{\alpha} (I)$ spaces are equal.
\end{thm}

We refer to (\cite{er} Theorem 2.13) for a proof. In fact, when $\alpha \not= n-1/2$, $n\in \N$ we have a stronger result as
the $J_{L,0}^{\alpha} (I)$, $J_{R,0}^{\alpha} (I)$ and $H_0^{\alpha} (I)$ spaces have equivalent semi-norms and norms.

\section{The extension problem}

As we want to deal with dynamical systems exhibiting the arrow of time, we need to consider the operator $\adta$ and $\tdbb$, in order to keep track of the
past and future of the dynamics. The fact that we consider $\alpha\not= \beta$ is only here for convenience. This can be used to take into account a
different quantity of {\it information} from the past and the future.\\

Let $\adta$ and $\tdbb$ be given. We look for an operator ${\cal D}^{\alpha ,\beta}$ of the form
\begin{equation}
{\cal D}^{\alpha ,\beta} =\M (\adta ,\tdbb ),
\end{equation}
where $\M :\R^2 \rightarrow \C$ is a mapping which does not depends on $(\alpha,\beta )$, satisfying the following general principles:
\vskip 3mm
\begin{itemize}
\item i) Gluing property: If $x(t)\in C^1$ then when $\alpha=\beta=m$, $m\in \N^*$, ${\cal D}^{m,m} x(t)=\di {d^m x\over dt^m}$.

\item ii) $M$ is a $\R$-linear mapping.

\item iii) Reconstruction: The mapping $M$ is {\it invertible}.
\end{itemize}
\vskip 3mm
Condition i) is fundamental in the embedding framework. It follows that all or constructions can be seen as a continuous {\it two-parameters deformation}
of the corresponding classical one\footnote{Condition i) is not the usual condition underlying the stochastic or quantum embedding theories.
In general, we have an injective mapping $\iota$ from the set of differentiable functions $C^1$ in a bigger functional space $\E$ such that the operator
$\cal D$ that we define on $\E$ reduces to the classical derivative on $\iota (C^1 )$ meaning that for $x\in C^1$, ${\cal D} (\iota (x))=\iota (x'(t))$
where $x'(t)=dx/dt$. As a consequence, we have a true embedding in this case, meaning that the embed theory already contain the classical one via
the mapping $\iota$. Here, the classical theory is not contained in the embedded theory but can be recovered by a continuous two-parameters deformation.}. This can be of importance
dealing with the {\it fractional quantization} problem in classical mechanics.\\

Condition ii) does not have a particular meaning. This is only the {\it simplest dependence} of the operator ${\cal D}$ with respect to
$\adta$ and $\tdbb$.\\

Condition iii) is important. It means that the data of ${\cal D}^{\alpha ,\beta}$ on a given function $x$ at point $t$ allows us to recover the left and right
$\RL$ derivatives of $x$ at $t$, so information about $x$ in a neighborhood of $x(t)$.

\begin{lem}
Operators satisfying conditions i), ii) and iii) are of the form
\begin{equation}
\label{contrainte1}
{\cal D}^{\alpha ,\beta} =[p \, \adta + (p-1) \, \tdbb ] +i q \, [\adta +\tdbb ] ,
\end{equation}
where $p$, $q\in \R$ and $q\not= 0$.
\end{lem}

\begin{proof}
By ii), we denote $M(x,y)=px+qy +i(rx+sy)$, with $p,q,r,s\in \R$. By i), we must have with $y=-x$ corresponding to the operator a choice
of operators $d/dt$, $-d/dt$
\begin{equation}
p-q=1,\ \ r-s=0 .
\end{equation}
We then already have an operator of the form (\ref{contrainte1}). The reconstruction assumption only impose that $q\not= 0$ in (\ref{contrainte1}).
\end{proof}

A more rigid form for these operators is obtained imposing a {\it symmetry} condition.
\vskip 3mm
\begin{itemize}
\item iv) Let $x \in C^0$ be a real valued function possessing left and right classical derivatives at point $t$, denoted by
$\di {d_+ x\over dt }$ and $\di {d_- x \over dt}$ respectively. If $\di {d_+ x\over dt }=-\di {d_- x \over dt}$, then we impose that
\begin{equation}
{\cal D}_{1,1} x(t)=i \di {d_+ x\over dt} .
\end{equation}
\end{itemize}
\vskip 3mm
Condition iv) must be seen as the non-differentiable pendant of condition i). Indeed, condition i) can be rephrased as follows: if $x\in C^0$ is such that
$d^+x$ and $d^- x$ exist and satisfy $d^+ x =d^- x$ then ${\cal D}^{1,1} x =d^+ x$. Condition iv) is then equivalent to the commutativity of the following
diagram, where $\R^2$ is seen as the $\R$-vector space associate to $\C$:
\begin{equation}
\left .
\begin{array}{lll}
\C & \stackrel{\tau}{\longrightarrow} & \C \\
x-ix & \longmapsto & x+ix \\
\M \downarrow & & \downarrow \M \\
\C & \stackrel{\tau}{\longrightarrow} & \C ,\\
a & \longmapsto & ia
\end{array}
\right .
\end{equation}
where $\tau :\C \rightarrow \C$ is defined by $\tau(z)=iz$, $z\in \C$ and we have used the fact that $M(x,-x)=0$ following condition i).

\begin{lem}
\label{extension}
The unique operator satisfying condition i), ii), iii) and iv) is given by
\begin{equation}
\label{contrainte2}
{\cal D}_{\alpha ,\beta} =\di {1\over 2} [\adta -\tdbb ] +i \di {1\over 2} [\adta +\tdbb ] .
\end{equation}
\end{lem}

\begin{proof}
By iv), we must have $p+(p-1)=0$ and $2q=1$, so that $p=q=1/2$.
\end{proof}

\section{The fractional operator of order $(\alpha ,\beta )$, $\alpha >0$, $\beta >0$}

Lemma \ref{extension} leads us to the following definition of a fractional operator of order $(\alpha ,\beta )$:

\begin{defi}
For all $a ,b \in \R$, $a<b$, the fractional operator of order $(\alpha ,\beta)$, $\alpha >0$, $\beta >0$, denoted by $\D$, is defined by
\begin{equation}
\D =\di {1\over 2}\left [ \adta -\tdbb \right ] +i\mu \di {1\over 2} \left [ \adta +\tdbb \right ] ,
\end{equation}
where $\mu \in \C$.
\end{defi}

When $\alpha=\beta=1$, we obtain $\D =d/dt$.\\

The free parameter $\mu$ can be used to reduce the operator $\D$ to some special cases of importance. Let us denoted by $x(t)$ a given
real valued function.\\

- For $\mu=-i$, we have $\D =\adta$ then dealing with an operator using the {\it future state} denoted by ${\cal F}_t (x)$ of the underlying
function, {\it i.e.} ${\cal F}_t (x)=\{ x(s),\ s\in [a,t[\}$.\\

- For $\mu=i$, we obtain $\D=-\tdbb$ then dealing with en operator using the {\it past state} denoted by ${\cal P}_t (x)$ of the underlying function,
{\it i.e.} ${\cal P}_t (x)=\{ x(s),\ s\in ]t,b]\}$.\\

As a consequence, our operator can be used to deal with problems using $\adta$, $\tdbb$, or both operators in a {\it unified} way, only particularizing
the value of $\mu$ at the end to recover the desired framework.\\

When $a=-\infty$ and $b=\infty$, we denote the associated operator ${\cal D}^{\alpha ,\beta}_{\mu }$ by $\FD^{\alpha ,\beta}_{\mu}$, {\it i.e.}
\begin{equation}
\FD^{\alpha ,\beta}_{\mu} =\di {1\over 2}\left [ \FD^{\alpha} -\FD^{\beta}_* \right ] +i\mu \di {1\over 2} \left [ \FD^{\alpha} +\FD^{\beta}_* \right ] ,
\end{equation}
where $\mu \in \C$.

\section{Product rules}

The classical product rule for Riemann-Liouville derivatives is for all $\alpha >0$
\begin{equation}
\label{prodruleini}
\di\int_a^b \adta f(t) g(t) dt =\int_a^b f(t) \tdba g(t) dt ,
\end{equation}
as long as $f(a)=f(b)=0$ or $g(a)=g(b)=0$.\\

This formula gives a strong connection between $\adta$ and $\tdba$ via a generalized integration by part. This relation is responsible for the emergence of
$\tdba$ in problems of fractional calculus of variations only dealing with $\adta$. See section \ref{fraceulerlagrange} for more details. This result also
justifies our approach to the construction of a fractional operator which put on the same level the left and right $\RL$ derivatives.\\

As a consequence, we obtain the following formula for our fractional operator:

\begin{lem}
For all $f,g \in \Ef$, we have
\begin{equation}
\label{prodrule}
\di\int_a^b \D f(t) g(t) dt =-\di\int_a^b f(t) \Di g(t) dt ,
\end{equation}
provide that $f(a)=f(b)=0$ or $g(a)=g(b)=0$.
\end{lem}

\begin{proof}
We have
\begin{equation}
\di\int_a^b \D f(t) g(t) dt =\di\int_a^b f (f) \left [ (\tdba -\adtb )+i\mu (\tdba +\adtb ) \right ] \left ( g(t)\right ) dt .
\end{equation}
Exchanging the role of $\alpha$ and $\beta$ in $(\tdba -\adtb )+i\mu (\tdba +\adtb )$, we obtain the operator
$(\tdbb -\adta )+i\mu (\tdbb +\adta )$ which can be written as
\begin{equation}
-\left [(\adta -\tdbb)-i\mu (\tdbb +\adta ) \right ] =-{\cal D}^{\alpha ,\beta}_{-\mu} .
\end{equation}
This concludes the proof.
\end{proof}

Here again, we see that it is convenient to keep the parameter $\mu$ free.

\newpage
\part{Fractional embedding of differential operators}
\label{part2}
\setcounter{section}{0}

\section{Fractional embedding of differential operators}

Let $d\in \N$ be a fixed integer and $a,b\in \R$, $a<b$ be given. We denote by $C([a,b])$ the set of continuous
functions $x:[a,b] \rightarrow \R^d$. Let $n\in \N$, we denote by $C^n ([a,b])$ the set of functions in $C([a,b])$ which are differentiable up to
order $n$.\\

Let $f :\R \times \C^d \rightarrow \C$ be a function, real valued on real arguments. We denote by $\F$ the corresponding operator acting on
functions $x$ and defined by
\begin{equation}
\F :
\left .
\begin{array}{lll}
C([a,b]) & \longrightarrow &  C ([a,b]) \\
x & \longmapsto & f(\bullet ,x(\bullet ) ) ,
\end{array}
\right .
\end{equation}
where $f(\bullet ,x(\bullet ))$ is the function defined by
\begin{equation}
f(\bullet ,x(\bullet ) ) :
\left .
\begin{array}{lll}
[a,b] & \longrightarrow & \C ,\\
t & \longmapsto & f(t, x(t)) .
\end{array}
\right .
\end{equation}
Let $\mathbf{f}=\{ f_i \}_{i=0,\dots ,n}$ be a finite family of functions, $f_i :\R \times \C^d \rightarrow \C$, and $F_i$, $i=1,\dots ,n$ the
corresponding family of operators. We denote by $\Op_{\mathbf{f}}$ the differential operator defined by
\begin{equation}
\label{formop}
\Op_{\mathbf{f}}^{\mathbf{g}} =\di\sum_{i=0}^n F_i \cdot \di {d^i \over dt^i } G_i ,
\end{equation}
where $\cdot$ is the standard product of operators, {\it i.e.} if $A$ and $B$ are two operators, we denote by $A\cdot B$ the operator defined by
$(A\cdot B)(x)=A(x) B(x)$ and $\circ$ the usual composition, {\it i.e.} $(A\circ B ) (x)=A(B(x))$, with the convention that $\di \left (
\di {d\over dt} \right )^0 =\Id$, where $\Id$ denotes the identity mapping on $\C$.

\begin{defi}[Fractional embedding of operators]
\label{defembop}
Let $\mathbf{f}=\{ f_i \}_{i=0,\dots ,n}$ and $\mathbf{g}=\{ g_i \}_{i=0,\dots ,n}$ be finite families of functions, $f_i :\R \times \C^d \rightarrow \C$ and
$g_i :\R \times \C^d \rightarrow \C$ respectively, and $F_i$, $G_i$, $i=1,\dots ,n$ the
corresponding families of operators, and $\Op_{\mathbf{f}}^{\mathbf{g}}$ the associated differential operator.\\

The $(\alpha ,\beta )$-fractional embedding of $\Op_{\mathbf{f}}^{\mathbf{g}}$ written as (\ref{formop}), denoted by $\Emb (\Op_{\mathbf{f}}^{\mathbf{g}} )$
is defined by
\begin{equation}
\label{formulaembop}
\Emb (\Op_{\mathbf{f}}^{\mathbf{g}} ) =\di\sum_{i=0}^n F_i \cdot \di \left ( \D \right )^i \circ G_i  .
\end{equation}
\end{defi}

Note that the embedding procedure acts on operators of {\it a given form} and not on operators like abstract data, {\it i.e.} this is not a mapping on
the set of operators.\\

We can solve this {\it indeterminacy} using a formal representation of an operator.\\

Let $\mathbf{f}=\{ f_i \}_{i=0,\dots ,n}$ and $\mathbf{g}=\{ g_i \}_{i=0,\dots ,n}$ be finite families of functions, $f_i :\R \times \C^d \rightarrow \C$ and
$g_i :\R \times \C^d \rightarrow \C$ respectively, and $F_i$, $G_i$, $i=1,\dots ,n$ the
corresponding families of operators. We denote by ${}_{\otimes}\Op_{\mathbf{f}}^{\mathbf{g}}$ the operator acting on $\E \otimes C^n \otimes \E$ defined by
\begin{equation}
\label{formoptensor}
{}_{\otimes} \Op_{\mathbf{f}}^{\mathbf{g}} =\di\sum_{i=0}^n F_i \otimes \di {d^i \over dt^i } \otimes \di G_i ,
\end{equation}
where $\otimes$ is the standard tensor product.

We denote by ${\cal O}_{\otimes}$ the set of operators of the form (\ref{formoptensor}) and $\cal O$ the set of differential operators of
the form (\ref{formop}). We define a mapping $\pi$ from ${\cal O}_{\otimes}$ to $\cal O$ by
\begin{equation}
\pi ({}_{\otimes} \Op_{\mathbf{f}}^{\mathbf{g}} )=\mu (\di\sum_{i=0}^n F_i \otimes \di (\D )^i \otimes G_i )=\sum_{i=0}^n F_i \cdot \di (\D )^i
\circ G_i ,
\end{equation}
where $\mu$ is the projection $\mu :\E \otimes \E \otimes \E \rightarrow \E$, $\mu (x\otimes y\otimes z)=x\cdot (y \circ z)$.\\

A differential operator being given, its fractional embedding depends on its writing as an element of ${\cal O}_{\otimes}$.

\section{Fractional embedding of differential equations}

Let $k\in \N$ be a fixed integer. Let $\mathbf{f}=\{ f_i \}_{i=0,\dots ,n}$ and $\mathbf{g}=\{ g_i \}_{i=0,\dots ,n}$ be finite families of functions,
$f_i :\R \times \C^{kd} \rightarrow \C$ and
$g_i :\R \times \C^{kd} \rightarrow \C$ respectively, and $F_i$, $G_i$, $i=1,\dots ,n$ the corresponding families of operators. We denote by
$\Op_{\mathbf{f}}^{\mathbf{g}}$ the operator acting on $(C^n [a,b] )^k$ defined by
\begin{equation}
\label{opequation}
\Op_{\mathbf{f}}^{\mathbf{g}} =\di\sum_{i=0}^n F_i \cdot \di {d^i \over dt^i } \circ \di G_i ,
\end{equation}
The ordinary differential equation associated to $\Op_{\mathbf{f}}^{\mathbf{g}}$ is defined by
\begin{equation}
\label{asequation}
\Op_{\mathbf{f}}^{\mathbf{g}} (x,\di {dx\over dt},\dots ,\di {d^k x\over dt^k } )=0 ,\ \ x\in C^{n+k} ([a,b]) .
\end{equation}

We then define the fractional embedding of equation (\ref{asequation}) as follow:

\begin{defi}
\label{defiembequation}
The fractional embedding of equation (\ref{asequation}) of order $(\alpha ,\beta )$, $\alpha ,\beta >0$ is defined by
\begin{equation}
\label{fracasequation}
\Emb \left ( \Op_{\mathbf{f}}^{\mathbf{g}} \right ) (x,\D x,\dots ,\di \left (\D \right )^k x  )=0 ,\ \ x\in \Ef (n+k) .
\end{equation}
\end{defi}

Note that as long as the form of the operator is fixed the fractional embedding procedure associates a {\it unique} fractional differential
equation.

\section{About time-reversible dynamics}

The fractional embedding procedure associate a natural fractional counterpart to a given ordinary differential equation. In some case, the underlying
equation possesses specific properties which have a physical meaning. One of this property is the {\it time-reversible} character of the dynamics:\\

A dynamics on a space $U$ is {\it time-reversible} if there exists an invertible map $i$ of $U$ such that $i^2 =\Id$, {\it i.e.} $i$ is an involution,
and if we denote by $\phi_t$ the flow describing the dynamics we have
\begin{equation}
i \circ \phi_{-t}  =\phi_t \circ i ,
\end{equation}
meaning that if $x(t)$ is a solution then $i (x(-t))$ is also a solution of the underlying equation.\\

Time-Reversibility is closely related to a specific property of the classical derivative under time-reversal:
\begin{equation}
\di {d \over dt} (x(-t))=- \di {dx\over dt} (-t) ,
\end{equation}

We define a notion of reversibility directly on operators:

\begin{defi}
We denote by $\Revo$ the $\C$- linear operator defined by $\Revo (\adta )=\tdba$, $\Revo (\tdba )=\adta$.
\end{defi}

The action of $\Revo$ on $\D$ is non-trivial:

\begin{lem}
\label{revactiond}
We have $\Revo (\D )= -\Di$.
\end{lem}

We then have the following analogue of reversibility in the fractional setting:

\begin{defi}
Let $\Op_{\mathbf{f}}^{\mathbf{g}}$ be a differential operator of the form (\ref{opequation}) such that the dynamics of the associated differential equation
(\ref{asequation}) is time-reversible, {\it i.e.}
\begin{equation}
\Revo (\Op_{\mathbf{f}}^{\mathbf{g}} )=-\Op_{\mathbf{f}}^{\mathbf{g}} .
\end{equation}
The fractional embedding $\Emb$ is called reversible if
\begin{equation}
\Revo (\Emb \left (\Op_{\mathbf{f}}^{\mathbf{g}}\right )  )=-\Emb \left ( \Op_{\mathbf{f}}^{\mathbf{g}} \right ).
\end{equation}
\end{defi}

The main consequence of lemma \ref{revactiond} is that there exists a unique way to do a fractional embedding conserving the reversibility symmetry.

\begin{thm}
The reversibility symmetry is preserved by a fractional embedding if and only if $\alpha=\beta$ and $\mu=0$.
\end{thm}

\begin{proof}
The reversibility symmetry is preserved if and only if we always have $\Revo (\D )=-\D$. By lemma \ref{revactiond} this is only possible when
$\mu=0$ and $\alpha=\beta$.
\end{proof}

In what follows we denote by $\Rev$ the fractional embedding ${}_a^{\alpha} \mbox{\rm Emb}_b^{\alpha} (0)$.

\newpage
\part{Fractional embedding of Lagrangian systems}
\label{part3}
\setcounter{section}{0}

In this section, we derive the fractional embedding of a particular class of ordinary differential equations called {\it Euler-Lagrange equations}
which governs the dynamics of {\it Lagrangian systems}.

\section{Reminder about Lagrangian systems}

Lagrangian systems play a central role in dynamical systems and physics, in particular for classical mechanics. We refer to \cite{ar} for more details.

\begin{defi}
An admissible Lagrangian function $L$ is a function $\L :\R \times \R^d \times \C^d \mapsto \C$ such that $\L (t,x,v)$ is holomorphic with respect to $v$,
differentiable with respect to $x$ and real when $v\in \R$.
\end{defi}

A Lagrangian function defines a {\it functional} on $C^1 (a,b)$, denoted by
\begin{equation}
\Ll_{a,b} : C^1 (a,b) \rightarrow \R , \ \ \ x\in C^1 (a,b) \longmapsto \int_a^b \L (s, x(s),\di {dx\over dt} (s) ) \ds ,\ \ a,b\in \R .
\end{equation}

The classical {\it calculus of variations} analyzes the behavior of $\Ll$ under small perturbations of the initial function $x$. The
main ingredient is a notion of differentiable functional and extremals.

\begin{defi}[Space of variations]
We denote by $\vari$ the set of functions in $C^1 (a,b)$ such that $h(a)=h(b)=0$.
\end{defi}

A functional $\Ll$ is {\it differentiable} at point $x\in C^1 (a,b)$ if and only if
\begin{equation}
\Ll (x+\epsilon h )-L(x)=\epsilon d\Ll (x,h) +o(\epsilon ) ,
\end{equation}
for $\epsilon >0$ and all $h\in \vari$.\\

Using the notion of differentiability for functionals one is lead to consider {\it extremum} of a given Lagrangian functional.

\begin{defi}
\label{classicalextremal}
An extremal for the functional $\Ll$ is a function $x\in C^1 (a,b)$ such that $d\Ll (x,h)=0$ for all $h\in \var$.
\end{defi}

Extremals of the functional $\Ll$ can be characterized by an ordinary differential equation of order $2$, called the Euler-Lagrange equation.

\begin{thm}
The extremals of $\Ll$ coincide  with the solutions of the Euler-Lagrange equation denoted by (EL) and defined by
\begin{equation}
\label{eulerlagrange}
\di {d\over dt} \left [ \di {\partial \L \over \partial v} \left ( t,x(t), \di {dx\over dt} (t)\right )\right ] =
\di {\partial L\over \partial x} \left ( t,x(t),\di {dx\over dt} (t)\right ) .
\end{equation}
\end{thm}

This equation can be seen as the action of the differential operator
\begin{equation}
\label{oplagrangian}
\mbox{\rm O}_{(EL)} =\di {d\over dt} \circ \di {\partial \L \over \partial v}  -\di {\partial L\over \partial x}
\end{equation}
on the couple $(x(t) , \di {dx\over dt} (t))$.\\

\section{Fractional Euler-Lagrange equation}

The fractional embedding procedure allows us to define a natural extension of the classical Euler-Lagrange equation in the fractional context. The
main result of this section is:

\begin{thm}
\label{defieulerlagrange}
Let $L$ be an admissible Lagrangian function. The $\Emb$-fractional Euler-Lagrange equation associated to $L$ is given by
$$
\di \D \left [ \di {\partial \L \over \partial v} \left ( t,x(t), \D x(t)\right )\right ] =
\di {\partial L\over \partial x} \left ( t,x(t),\D x(t)\right ) .
\eqno{(\FEL_{\alpha ,\beta}^{\mu} )}
$$
\end{thm}

In what follows, we will simply speak about the fractional Euler-Lagrange equation when there is no confusion on the underlying embedding procedure.\\

The proof is based on the following lemma:

\begin{lem}
\label{lememboplagrangian}
Let $\L$ be an admissible Lagrangian function. The fractional embedding of the Euler-Lagrange differential operator $\mbox{\rm O}_{(EL)}$ is given by
\begin{equation}
\label{emboplagrangian}
\Emb (\mbox{\rm O}_{(EL)} ) =\di \D \circ \di {\partial \L \over \partial v} -\di {\partial L\over \partial x} .
\end{equation}
\end{lem}

\begin{proof}
The operator (\ref{oplagrangian}) is first considered as acting on $(C^1 ([a,b])^2$, {\it i.e.} for all $(x(t),y(t)) \in C^1 [a,b] \times C^1[a,b]$
we have
\begin{equation}
\mbox{\rm O}_{(EL)} (x(t),y(t)) =\di {d\over dt} \left ( \di {\partial \L \over \partial v} (t,x(t),y(t)) \right )-\di {\partial L\over \partial x}
(t,x(t),y(t)) .
\end{equation}
This operator is of the form $\mbox{\rm O}_{\mathbf{f}}^{\mathbf{g}}$ with
\begin{equation}
\mathbf{f} = \left ( \mbox{\rm\bf 1} , \di {\partial L\over \partial x} \right ) ,
\end{equation}
and
\begin{equation}
\mathbf{g} =\left ( -\di {\partial \L \over \partial v} ,\mbox{\rm\bf 1} \right ) ,
\end{equation}
where $\mbox{\rm\bf 1} :\R \times \C^2 \rightarrow \C$ is the {\it constant} function $\mbox{\rm\bf 1} (t,x,y)=1$. As a consequence, $\mbox{\rm O}_{(EL)}$
is given by
\begin{equation}
\mbox{\rm O}_{(EL)}= \mbox{\rm\bf 1} \cdot \di {d\over dt} \circ \di {\partial L\over \partial v }
-\di {\partial L\over \partial x} \cdot \di \Id \circ \mbox{\rm\bf 1} ,
\end{equation}
with the convention that $\di \left ( \di {d\over dt} \right )^0 =\Id$. We then obtain equation (\ref{emboplagrangian}) using definition \ref{defembop}.
\end{proof}

We can now conclude the proof of theorem \ref{defieulerlagrange} using definition \ref{defiembequation}. The fractional embedding of equation
(\ref{eulerlagrange}) is given by
\begin{equation}
\Emb \left (\mbox{\rm O}_{(EL)} \right ) \left ( x ,\D x \right )=0 ,
\end{equation}
which reduces to equation $(\FEL_{\alpha ,\beta}^{\mu} )$ thanks to lemma \ref{lememboplagrangian}.

\section{The coherence problem}
\label{sectioncoherence}

The fractional embedding procedure allows us to define a natural fractional analogue of the Euler-Lagrange equation. This result is satisfying
because the procedure is fixed. However, Lagrangian systems possess very special features. In particular, the classical Euler-Lagrange
equation can be obtained using a {\it variational principle}, called the {\it least-action principle} and denoted $\lap$. The least
action principle asserts that the Euler-Lagrange equation characterizes the extremals of a given functional associated to the Lagrangian. We
then are lead to the following problem:
\vskip 3mm
\begin{itemize}
\item i) Develop a calculus of variation on fractional functionals.

\item ii) State the corresponding {\it fractional least-action principle}, in particular explicit the associated {\it fractional Euler-Lagrange equation}
denoted by $\FEL_{flap}$.

\item iii) Compare the result with the embedded Euler-Lagrange equation ($\FEL_{\alpha ,\beta}^{\mu}$)
\end{itemize}
\vskip 3mm
An embedding procedure is called {\it coherent} when the two Euler-Lagrange equations are the same, {\it i.e.} if
\begin{equation}
\FEL_{flap} =\Emb (EL) ,
\end{equation}
assuming that $\FEL_{flap}$ is obtained from the embedding of the classical functional using the {\it same} embedding procedure.\\

As we will see, an embedding procedure is not always coherent. Although we obtain in general equations of the same form, we usually have some
{\it torsion} between the embedding of the functional and the Euler-Lagrange equation which cancel only in particular cases.\\

The fractional calculus of variations is developed in $\S$.\ref{fraccalculvar} as well as the corresponding fractional least-action principle. The coherence
of fractional embedding procedures is discussed in $\S$.\ref{coherence}.

\newpage
\part{Fractional calculus of variations}
\setcounter{section}{0}
\label{fraccalculvar}

This section is devoted to the {\it fractional calculus of variations} using our fractional operator. The functional is obtained under the fractional
embedding procedure. We refer to the work of O.P. Agrawal \cite{agra} for related results.

\section{Fractional functional}

Let $\L$ be an admissible Lagrangian function on $\R\times \R^d \times \C^d$, $d\geq 1$, and $\Ll$ the associated functional.
Using the Fractional embedding procedure $\Emb$, we define a natural {\it Fractional functional} associated to $L$.\\

We denote by $\Ef$ the set of functions $x$ such that $\adta x$ and $\tdbb x$ are defined.

\begin{defi}
The Fractional functional associated to $L$ is defined by
\begin{equation}
\Ll_{a,b}^{\alpha ,\beta} :\Ef \rightarrow \R , \ \ \ x\in \Ef \longmapsto
\int_a^b \L (s, x(s),\D x (s) ) \ds ,\ \ a,b\in \R ;
\end{equation}
\end{defi}

The extension property implies that $\Ll_{a,b}^{\alpha ,\beta}$ reduce to the classical functional $\Ll_{a,b}$ when $\alpha=\beta =1$.

\section{Space of variations and extremals}

Let us denote by $\Ef (0,0)$ the set of curves $h\in \Ef$ satisfying $h(a)=h(b)=0$ and ${}_a E^{\alpha}_b :={}_a^{\alpha} E^{\alpha}_b$. We denote by
$\Var^{\alpha} (0,0)$ the set defined to be
\begin{equation}
\Var^{\alpha} (a,b)= \{ h\in {}_a E^{\alpha}_b ,\ h(a)=h(b)=0\ \ \mbox{\rm and}\ \ \adta h =\tdba h \ \} .
\end{equation}

We denote by $\P$ the set $\Ef (0,0)$ or $\Var^{\alpha} (a,b)$.

\begin{defi}
Let $x$ be a given curve. A $\P$-variation of $x$ is a one-parameter $\epsilon \in \R$ family of curves of the form
\begin{equation}
y_{\epsilon} =x+\epsilon h ,\ \ h\in \P .
\end{equation}
\end{defi}

A notion of differentiability can now be defined for fractional functionals. In the following, we write $\Ll_{a,b}$ indifferently for $\Ll_{a,b}^{\alpha ,\beta}$ when
$\P =\Ef (0,0)$ and $\Ll_{a,b}^{\alpha,\alpha}$ when $\P =\Var^{\alpha} (a,b)$.

\begin{defi}
Let $L$ be an admissible Lagrangian function and $\Ll_{a,b}$ the associated fractional functional. The functional $\Ll_{a,b}$
is called $\P$-differentiable at $x$ if
\begin{equation}
\Ll_{a,b} (x+\epsilon h ) - \Ll_{a,b} (x)= \epsilon d\Ll_{a,b} (x,h) +o(\epsilon ),
\end{equation}
for all $h\in \P$, $\epsilon >0$, where $d\Ll_{a,b} (x,h)$ is a linear functional of $h$.
\end{defi}

The linear functional $d\Ll_{a,b}(x,h)$ is called the $\P$-{\it differential} of the fractional functional $\Ll_{a,b}$ at point $x$.\\

An {\it extremal} for $\Ll_{a,b}$ is then defined by:

\begin{defi}
A $\P$-extremal for the functional $\Ll_{a,b}$ is a function $x$ such that $d\Ll_{a,b}  (x,h)=0$ for all
$h\in \P$.
\end{defi}

The following lemma gives the explicit form of the differential of a fractional functional:

\begin{thm}
\label{tdifferential}
Let $L$ be an admissible Lagrangian function and $\Ll_{a,b}^{\alpha ,\beta}$ the associated fractional functional. The functional
$\Ll_{a,b}^{\alpha ,\beta}$ is differentiable at any $x\in \Ef (x_a ,x_b )$ and for all $h\in \Ef (0,0)$ the differential is given by
\begin{equation}
\label{differential}
d\Ll_{a,b}^{\alpha ,\beta} (x,h) = \int_a^b \left [- \Di \left [
\di {\partial \L\over \partial v} \left ( t,x(t),\D x (t) \right ) \right ] +\di {\partial \L \over \partial x} (t,x(t),\D x(t))
\right ] h(t) dt .
\end{equation}
\end{thm}

\begin{proof}
As the left and right $\RL$ derivatives are linear operators we have
\begin{equation}
\D (x+\epsilon h)=\ D x +\epsilon \,\D h.
\end{equation}
As a consequence, we obtain
\begin{equation}
\Ll_{a,b}^{\alpha ,\beta} (x+ \epsilon h) =\int_a^b \L (s, x(s)+\epsilon h(s),\D x (s) +\epsilon\, \D h (s)) \ds
\end{equation}
which implies, doing a Taylor expansion of $\L (s, x(s)+\epsilon h(s),\D x (s) +\epsilon\, \D h (s))$ in $\epsilon$ around $0$
\begin{equation}
\label{forminter}
\left .
\begin{array}{lll}
\Ll_{a,b}^{\alpha ,\beta} (x+ \epsilon h) & = & \epsilon \di\int_a^b \di \left [
\di {\partial \L \over \partial x} (s, x(s),\D x (s)) h(s) +\di {\partial \L \over \partial v} (s, x(s),\D x (s)) \D h(s) \right ] \ds \\
 & & + \Ll_{a,b}^{\alpha ,\beta} (x) +o(\epsilon ).
\end{array}
\right .
\end{equation}
Using the product rule (\ref{prodrule}) we obtain
\begin{equation}
\int_a^b \di {\partial \L \over \partial v} (s, x(s),\D x (s)) \D h(s) \ds =-\int_a^b \di \Di {\partial \L \over \partial v} (s, x(s),\D x (s)) h(s) \ds .
\end{equation}
Replacing this expression in (\ref{forminter}), we deduce formula (\ref{differential}).
\end{proof}

\section{The fractional Euler-Lagrange equation}
\label{fraceulerlagrange}

We obtain the following analogue of the {\it least-action principle} in classical Lagrangian mechanics:

\begin{thm}[Fractional least-action principle]
\label{teulerlagrange}
Let $\Ll [x]$ be a functional of the form
\begin{equation}
\Ll [x]=\int_a^b \L (s, x(s) ,\D x(s) ) ds
\end{equation}
defined on $\Ef (x_a ,x_b)$.

A necessary and sufficient condition for a given function $x\in \Ef$ to be a $\Ef$-extremal for $\Ll [x]$ with fixed end points $x(a)=x_a$, $x(b)=x_b$,
is that it satisfies the fractional Euler-Lagrange equation (FEL):
\begin{equation}
\label{fraceulerlagrange2}
\Di \left [
\di {\partial \L\over \partial v} \left ( t,x(t),\D x (t) \right ) \right ] =\di {\partial \L \over \partial x} (t,x(t),\D x(t)) .
\end{equation}
\end{thm}

Note that this equation is different from the one obtained via  the fractional embedding procedure.

\begin{proof}
Using the classical Du Bois Reymond lemma (\cite{av},p.108) and theorem \ref{tdifferential} we obtain (\ref{fraceulerlagrange2}).
\end{proof}

The weak analogue using the space of variation $\Var^{\alpha} (a,b)$ is given by:

\begin{thm}[Weak fractional least-action principle]
\label{wteulerlagrange}
Let $\Ll [x]$ be a functional of the form
\begin{equation}
\Ll [x]=\int_a^b \L (s, x(s) ,{\cal D}^{\alpha}_{\mu} x(s) ) ds
\end{equation}
defined on ${}_a E^{\alpha}_b (x_a ,x_b)$.

A necessary and sufficient condition for a given function $x\in {}_a E^{\alpha}_b$ to be a $\Var^{\alpha} (a,b)$-extremal for $\Ll [x]$ with fixed end points
$x(a)=x_a$, $x(b)=x_b$, is that it satisfies the fractional Euler-Lagrange equation $\FEL^{\alpha ,\alpha}_{\mu}$
\end{thm}

We denote $\FEL^{\alpha}_{\mu}$ for $\FEL^{\alpha ,\alpha}_{\mu}$ in the following.

\section{Coherence}
\label{coherence}

The coherence problem can now be studied in details. We have the following theorem, which is only a rewriting of theorem \ref{teulerlagrange} and
definition \ref{defieulerlagrange}:

\begin{thm}
\label{resume}
Let $L$ be an admissible Lagrangian function, $a,b\in \R$, $a<b$, $\alpha ,\beta >0$, then the following diagram commutes:
\begin{equation}
\left .
\begin{array}{lll}
\L (t, x(t), dx/dt ) & \stackrel{\Emb}{\longrightarrow} & \L (t,x(t) ,\D x )\\
\lap \downarrow & & \downarrow \flap\\
\EL & \stackrel{\Embi}{\longrightarrow} & \FEL_{\beta ,\alpha }^{-\mu} .
\end{array}
\right .
\end{equation}
\end{thm}

Theorem \ref{resume} is not a {\it coherence} result in the spirit of (\cite{cr2},\cite{cr3}) or \cite{cd3}. Indeed, the embedding procedure changes
between the Euler-Lagrange equation and the functional.\\

There exist at least two ways to restore coherence of the embedding procedure.\\

The first one is to restrict the set of variations we are looking for in the least-action principle. Using the set $\Var^{\alpha} (a,b)$ we obtain a coherent
procedure via the weak least action principle.

\begin{thm}[Weak coherence]
Let $L$ be an admissible Lagrangian function, $a,b\in \R$, $a<b$, $\alpha>0$, then the following diagram commutes:
\begin{equation}
\left .
\begin{array}{lll}
\L (t, x(t), dx/dt ) & \stackrel{{}_a \mbox{\rm Emb}^{\alpha}_b}{\longrightarrow} & \L (t,x(t) ,{\cal D}^{\alpha}_{\mu} x )\\
\lap \downarrow & & \downarrow \mbox{\rm WLAP}\\
\EL & \stackrel{{}_a \mbox{\rm Emb}^{\alpha}_b}{\longrightarrow} & \FEL^{\alpha }_{\mu} ,
\end{array}
\right .
\end{equation}
where $\mbox{\rm WLAP}$ stands for the weak least-action principle and ${}_a \mbox{\rm Emb}^{\alpha}_b :={}_a^{\alpha} \mbox{\rm Emb}^{\alpha}_b$.
\end{thm}

The second way is to use a specific fractional operator with a particular symmetry property. Using the reversible embedding and the
operator $\mbox{\rm d}^{\alpha}$ which corresponds to the operator ${\cal D}^{\alpha ,\alpha}_{0}$, {\it i.e.}
\begin{equation}
\mbox{\rm d}^{\alpha} := \di {\adta +\tdba \over 2} ,
\end{equation}
we obtain the following coherence result:

\begin{thm}[Reversible coherence]
Let $L$ be an admissible Lagrangian function, $a,b\in \R$, $a<b$, $\alpha>0$, then the following diagram commutes:
\begin{equation}
\left .
\begin{array}{lll}
\L (t, x(t), dx/dt ) & \stackrel{\Rev}{\longrightarrow} & \L (t,x(t) ,\mbox{\rm d}^{\alpha} x )\\
\lap \downarrow & & \downarrow \flap\\
\EL & \stackrel{\Rev}{\longrightarrow} & \FEL^{\alpha } .
\end{array}
\right .
\end{equation}
where $\Rev$ is the $\alpha$-reversible embedding and $\FEL^{\alpha}$ denotes the equation $\FEL^{\alpha ,\alpha}_0$.
\end{thm}

It is not clear for the moment to know what is the interest of the fractional reversible embedding procedure for applications. At least it can be
considered as a good candidate to a {\it fractional quantization} procedure.

\newpage
${}$
\newpage
\part{Symmetries and the Fractional Noether theorem}
\label{part5}
\setcounter{section}{0}

In classical mechanics Noether's theorem gives a strong connection between group of symmetries of the Lagrangian and conservation laws. We refer
in particular to Arnold's presentation of Noether's theorem for Lagrangian systems (\cite{ar} p.88). For a historical point of view on different
generalization of Noether's theorem we refer to \cite{sch}.

\section{Invariance of fractional functionals}

Group of symmetries are classical in mechanics. A natural way to deal with symmetries is to look for one-parameter group of diffeomorphisms. We first
define the action of a diffeomorphism on a couple $(x(t), \D x(t))$.

\begin{defi}
Let $\phi :\R^d \rightarrow \R^d$ be a diffeomorphism. The fractional linear tangent map associated to $\phi$, denoted by $\phi_* $ is defined by
\begin{equation}
\phi_* (x, \D x) =(\phi(x), \D \phi(x)) .
\end{equation}
\end{defi}

\begin{defi}
Let $\Phi =\{ \phi_s \} _{s\in \R}$ be a one-parameter family of diffeomorphisms. An admissible Lagrangian $\L$ is said to be invariant under the action of
$\Phi$ if
\begin{equation}
\L (x,\D x)=\L (\phi_s (x) ,\D \left (\phi_s (x) \right ) ) ,\ \ \forall s\in \R ,
\end{equation}
or equivalently
\begin{equation}
\L (x,\D x )=\L \left ( (\phi_s )_* (x,\D x) \right ) .
\end{equation}
\end{defi}

The relation with {\it infinitesimal transformations} used in \cite{gasto}, is obtained using a Taylor expansion of $y_t (s)=\phi_s (x(t))$ in a neighborhood of $0$,
leading to
\begin{equation}
y_t (s)=y_t (0)+ s. \di {dy_t \over ds} (0) +o(s) .
\end{equation}
As $\phi_0 =\Id$ is the identity map, we obtain denoting $\di {dy_t \over ds} (0) =\xi (t,x)$ an infinitesimal transformation given by
\begin{equation}
x(t) \mapsto x(t)+s \xi (t,x(t)) +o(s) .
\end{equation}

A natural problem in the embedding framework is the following: assume that $\L$ is invariant under a one-parameter group of diffeomorphisms. What can
we said about the fractional embedded Lagrangian ?

\section{The fractional Noether theorem}

In this paragraph, we derive a fractional version of the Noether theorem. Similar results have been obtained by G. Frederico and D. Torres \cite{gasto} in
a different setting.

\begin{thm}[The fractional Noether theorem]
Let $\L$ be an admissible Lagrangian invariant under a one-parameter group of diffeomorphism $\Phi =\{ \phi_s \}_{s\in \R}$. Then, we have
\begin{equation}
\label{noetherform}
\Di \left ( \di {\partial \L \over \partial v} \right ) \di {d y_t \over d s} +
\di {\partial \L \over \partial v} \di \D \left ( \di {d y_t \over d s} \right ) \mid_{s=0} =0 ,
\end{equation}
where $y_t (s)=\phi_s (x(t))$, along all solutions of the fractional Euler-Lagrange equation.
\end{thm}

Equation (\ref{noetherform}) is not the usual way to state the Noether theorem in classical mechanics. This is mainly due to the fact that the
fractional differential operator $\D$ inherits from the bad properties of the underlying left and right Riemann-Liouville derivatives concerning
the Leibniz rule.

\begin{proof}
Let $y_t (s)=\phi_s (x(t))$, the $\Llf$ invariance is equivalent to
\begin{equation}
\di {d\over ds} \left [ \L (t,y_t (s) ,\D y_t (s)) \right ] =0.
\end{equation}
The usual chain rule for the classical derivative implies
\begin{equation}
\di {\partial \L \over \partial x} (t,y_t (s) ,\D y_t (s)) \cdot \di {dy_t \over ds} +\di {\partial \L \over \partial v} (t,y_t (s) ,\D y_t (s))
\cdot \di {d\over ds} \left [ \D y_t (s) \right ] =0 .
\end{equation}
As $\D$ acts on the variable $t$ and $d/ds$ on the variable $s$, and $\di {dy_t (s) \over ds} \in \Ef$, we deduce that
\begin{equation}
\di {d\over ds} \left [ \D y_t (s) \right ] =\D \left [ \di {dy_t \over ds} \right ] .
\end{equation}
As a consequence, we obtain
\begin{equation}
\di {\partial \L \over \partial x} (t,y_t (s) ,\D y_t (s)) \cdot \di {dy_t \over ds} +\di {\partial \L \over \partial v} (t,y_t (s) ,\D y_t (s))
\cdot \di \D \left [ \di {d y_t \over ds} \right ] =0 .
\end{equation}
As $x(t)$ is an extremal for $\L$, we have
\begin{equation}
\Di \left ( \di {\partial \L \over \partial v} \right ) =\di {\partial \L \over \partial x} .
\end{equation}
This concludes the proof.
\end{proof}

\section{Toward fractional integrability and conservation laws}

Of course, when $\alpha=\beta=1$, equation (\ref{noetherform}) reduces to
\begin{equation}
\di {d\over dt} \left ( \di {\partial \L \over \partial v} \right ) \di {d y_t \over d s} +
\di {\partial \L \over \partial v} \di {d\over dt} \left ( \di {d y_t \over d s} \right ) \mid_{s=0} =0 ,
\end{equation}
which is equivalent to
\begin{equation}
\di {d\over dt} \left ( \di {\partial \L \over \partial v} \di {d y_t \over ds} (0) \right )  =0 ,
\end{equation}
using the Leibniz rule.\\

In that case, the Noether theorem can then be stated in term of {\it conservation laws} (first integrals), {\it i.e.} saying that the
quantity
\begin{equation}
C(t,x(t) , \dot{x} (t)) =\di {\partial \L \over \partial v} (t,x(t),\dot{x} (t)) \cdot \di {d y_t \over ds} (0) ,
\end{equation}
is {\it constant} along all the solutions of the Euler-Lagrange equation.\\

First integrals play a fundamental role in classical Lagrangian dynamics
as they are related to the classical problem of {\it integrability} (Liouville's integrability), {\it i.e.} the fact that a Lagrangian systems with
sufficiently many first integrals can be integrated by {\it quadratures} (see \cite{ar} Chapter 10).\\

It is not clear regarding to the fractional Noether's theorem to define the analogue of conservation laws and in particular the fractional analogue
of integrability for Lagrangian systems.

\newpage
${}$
\newpage
\part{Fractional Hamiltonian systems}
\label{part6}
\setcounter{section}{0}

The classical passage from the {\it configuration space}, {\it i.e.} space of positions to the {\it phase space}, {\it i.e.} the space of position and
velocities, is achieved using the Hamiltonian formalism. In this section, we derive the fractional analogue of Hamiltonian systems.

\section{Hamiltonian systems and Legendre property}

Let $L$ be a Lagrangian system. When the Lagrangian satisfies an analytic property, that we called the {\it Legendre property}, we can find a
remarkable symmetric formulation of the Euler-Lagrange equation called {\it Hamiltonian}. The Legendre property can be stated as follows:

\begin{defi}
Let $L$ be an admissible Lagrangian function. We say that $L$ satisfies the Legendre property if the mapping $v\mapsto \di {\partial L\over \partial v} (x,v)$
is invertible for all $x\in \R^d$.
\end{defi}

As a consequence, there exists a map $f:\R^d \times \C^d \rightarrow \C^d$ such that for any $p\in \C^d$ satisfying $p=\di {\partial L\over \partial v} (x,v)$
we have $v=f(x,p)$. The map $f$ is called the {\it Legendre transform}.\\

The Legendre transform connects the {\it momentum}
\begin{equation}
p(t)=\di {\partial L\over \partial v} (x(t),\di \dot{x} (t) ) ,
\end{equation}
of a given solution $x(t)$ of the Euler-Lagrange equation with its velocity $\dot{x} (t)$.\\

Using the Legendre transform, we can define a new function called the {\it Hamiltonian} associated to $L$.

\begin{defi}
Let $L$ be an admissible Lagrangian system satisfying the Legendre property, and $f$ the corresponding Legendre transform. The Hamiltonian function
associated to $L$, denoted by $H:\R^d \times \C^d \rightarrow \C$, is defined by
\begin{equation}
\label{ham}
H(x,p)=pf(x,p)-L(x,f(x,p)) .
\end{equation}
\end{defi}

The Hamiltonian system associated to $H$ is given as
\begin{equation}
\left .
\begin{array}{lll}
\di {dx\over dt} & = & \di {\partial H\over \partial p} ,\\
\di {dp\over dt} & = & - \di {\partial H\over \partial x} .
\end{array}
\right .
\end{equation}
The dynamics of the Hamiltonian system is equivalent to the dynamics governed by the Lagrangian.

\section{The fractional momentum and Hamiltonian}

The previous construction can be carried in the fractional context for our fractional Euler-Lagrange equation.

\begin{defi}
The fractional momentum of a given solution $x(t)$ of the fractional Euler-Lagrange equation is defined as
\begin{equation}
\p(t)=\di {\partial L\over \partial v} (x(t),\D x(t)) .
\end{equation}
\end{defi}

We remark that the fractional momentum is obtain from the classical one by the fractional embedding procedure. If $L$ satisfies the
Legendre property we obtain a fractional analogue of the relation between momentum and velocities:

\begin{lem}
Let $L$ be an admissible Lagrangian system satisfying the Legendre property and $f$ the associated Legendre transform. Let $x(t)$ be a solution of the
fractional Euler-Lagrange equation, and $\p (t)$ its fractional momentum. We have
\begin{equation}
\D x(t)=f(x(t),\p(t)) .
\end{equation}
\end{lem}

The fractional Euler-Lagrange equation is then equivalent to the following system of fractional differential equations:
\begin{equation}
\label{frachamiltonian1}
\left .
\begin{array}{lll}
\D x (t) & = & f (x(t), \p(t) ) ,\\
\Di \p (t) & = & \di {\partial L\over \partial x} (x(t) ,f(x(t),\p (t))) ,
\end{array}
\right .
\end{equation}
which is nothing else that
\begin{equation}
\label{frachamiltonian}
\left .
\begin{array}{lll}
\D x (t) & = & \di {\partial H\over \partial p}f (x(t), \p (t) ) ,\\
\Di \p (t) & = & -\di {\partial H\over \partial x} (x(t) ,\p (t)) ,
\end{array}
\right .
\end{equation}
with $H$ the Hamiltonian function (\ref{ham}).\\

We then introduce the following definition of a fractional Hamiltonian system:

\begin{defi}
Let $L$ be an admissible Lagrangian system satisfying the Legendre property. The fractional Hamiltonian system associated to $L$ is defined by
equation (\ref{frachamiltonian}).
\end{defi}

\section{Fractional Hamilton least-action principle}

Fractional Hamilton equations (\ref{frachamiltonian}) can also be obtained via a variational principle called {\it fractional Hamilton leat-action
principle}.

\begin{thm}
We denote by $H:\R^ d \times \C^d \rightarrow \C$ a Hamiltonian function and $F$ the associated functional defined by
\begin{equation}
\label{funchamilton}
F (x(t),\p (t))=\di\int_a^b \left  (\p (t) \D x(t) -H(x(t) ,\p (t)) \right ) dt .
\end{equation}
A couple $(x(t),\p (t))$ such that $x(a)=x_a$, $\p (a)=\p_a$, $x(b)=x_b$, $p(b)=\p_b$ is an extremal of (\ref{funchamilton}) if and only if it satisfies
the fractional Hamiltonian equations (\ref{frachamiltonian}).
\end{thm}

\begin{proof}
We denote by  $\mathbf{L}$ the function
\begin{equation}
\mathbf{L} (x,\p ,v,w)=\p v-H(x,\p ) .
\end{equation}
The functional $F$ can be seen as a functional associated to the Lagrangian $\mathbf{L}$. The fractional Euler-Lagrange equation associated to this
functional is
\begin{eqnarray}
\Di \left ( \di {\partial \mathbf{L} \over \partial v} (z(t)) \right ) =\di {\partial \mathbf{L} \over \partial x} (z(t)) ,\\
\Di \left ( \di {\partial \mathbf{L} \over \partial w} (z(t)) \right ) =\di {\partial \mathbf{L} \over \partial \p} (z(t)) ,\\
\end{eqnarray}
where $z(t)=(x(t),\p (t) ,\D x(t) ,\D \p (t))$.\\

As $\di {\partial \mathbf{L} \over \partial v} =\p$ and $\di {\partial \mathbf{L} \over \partial x} =-\di {\partial H\over \partial x}$, the first equation
gives
\begin{equation}
\Di \p (t) =-\di {\partial H\over \partial x} (x(t),\p (t)) .
\end{equation}
Moreover, as $\di {\partial \mathbf{L}\over \partial w} =0$ and $\di {\partial \mathbf{L} \over \partial p} =v-\di {\partial H\over \partial \p}$,
the second equation gives
\begin{equation}
0=\D x(t)-\di {\partial H\over \partial \p } (x(t),\p (t)) .
\end{equation}
This concludes the proof.
\end{proof}

As a consequence, the usual construction of Hamiltonian equations using the Legendre transform and fractional momentum is coherent with the fractional
Hamilton least-action principle.

\newpage
\part{Fractional Ostrogradski formalism and Nonconservative dynamical systems}
\label{part7}
\setcounter{section}{0}

This part is devoted to the fractional embedding of functionals depending on higher order derivatives, {\it i.e.} to a fractional analogue of the
Ostrogradski formalism. We derive the generalized fractional
Euler-Lagrange equation. As an application, we prove that some Nonconservative dynamical systems can be obtained via the fractional embedding of
well chosen higher order Lagrangian systems. In particular, we obtain a new approach to results obtain previously by F. Riewe (\cite{ri1} \cite{ri2}).

\section{Ostrogradski formalism and fractional embedding}

We first introduce some convenient definitions and notations.

\begin{defi}
Let $n\in \N^*$ and $d\in \N^*$, a Lagrangian function of order $n$ denoted by $\L$ is a function defined by
\begin{equation}
\left .
\begin{array}{lll}
\R \times \C^d \times \underbrace{\C^d \times \dots \times \C^d}_{n\ \mbox{\rm times}} & \longrightarrow & \C ,\\
(t,x,v_1 ,\dots ,v_n ) & \longmapsto & \L (t,x,v_1 ,\dots ,v_n ) ,
\end{array}
\right .
\end{equation}
for $x\in \C^d$, $v_i \in \C^d$, $i=1,\dots ,n$.
\end{defi}

A Lagrangian function of order $n$ possesses a natural Lagrangian functional denoted by $\Ll$ and defined by
\begin{equation}
\Ll : x(t) \longmapsto \int_a^b \L (t,x(t) , x^{(1)} (t) ,\dots ,x^{(n)} (t) ) dt ,
\end{equation}
where $x^{(n)} (t)$ denotes the $n$-th derivative of $x$ with respect to $t$, {\it i.e.}
\begin{equation}
x^{(n)} =\di {d^n x \over dt^n } =\di \left ( \di {d\over dt} \right )^n x ,
\end{equation}
where
\begin{equation}
\di \left ( \di {d\over dt} \right )^n =\underbrace{\di {d\over dt} \circ \dots \circ \di {d\over dt}}_{n\ \mbox{\rm times}} ,
\end{equation}
with $\circ$ the usual composition of operators.\\

The associated notion of extremals differs from definition \ref{classicalextremal} only on the boundary conditions, {\it i.e.} on the {\it space of variations}.

\begin{defi}[Space of variations of order $n$]
Let $n\in \N^*$, we denote by $\var^n $ the subset of functions $h\in C^n (a,b)$ satisfying
\begin{equation}
\left .
\begin{array}{l}
h(a)=h'(a)=\dots =h^{(n-1)} (a)=0 ,\\
h(b)=h'(b)=\dots =h^{(n-1)} (b)=0 .
\end{array}
\right .
\end{equation}
\end{defi}

We deduce the following generalized notion of extremals.

\begin{defi}
\label{genextremals}
Let $n\in \N^*$, $\L$ be a Lagrangian function of order $n$ and $\Ll$ the associated functional. The functional $\Ll$ is called differentiable at
$x \in C^n (a,b)$ if and only if
\begin{equation}
\Ll (x+\epsilon h ) -\Ll (x)=\epsilon d \Ll (x,h) +o(\epsilon ) ,
\end{equation}
for all $h\in \var^n$ where $d\Ll (x,h)$ is a linear functional of $h$.
\end{defi}

We then define extremals of the functional $\Ll$ as follows:

\begin{defi}
An extremal for the functional $\Ll$ is a function $x\in C^n (a,b)$ such that d$\Ll (x,h)=0$ for all $x\in \var^n$.
\end{defi}

Extremals of the functional $\Ll$ are characterized by a {\it generalized Euler-Lagrange} equation given by
$$
\di {\partial L\over \partial x} (t,x(t) ,x^{(1)} (t) ,\dots ,x^{(n)} (t)) + \di \sum_{i=1}^n (-1)^i \di \left ( \di  {d\over dt} \right ) ^i
\left [ {\partial L\over \partial v_i} (t,x(t),x^{(1)} (t) ,\dots ,x^{(n)} (t)) \right ] .
\eqno{(GEL)}
$$
We refer to (\cite{gf} Chap.2 $\S$.11 p.42) for a proof. We denote by OLAP the {\it Ostrogradski least-action principle}, {\it i.e.} the
derivation of (GEL) from $\Ll$.\\

The generalized Euler-Lagrange equation can be written using a differential operator denoted by $\Op_{\rm GEL}$ and defined by
\begin{equation}
\Op_{\rm GEL} =\di {\partial L\over \partial x } +\di\sum_{i=1}^n (-1)^i \left ( \di {d\over dt} \right )^i \circ \di {\partial L\over \partial v_i} .
\end{equation}
The (GEL) is then equivalent to
\begin{equation}
\label{geneulerlagrange}
\Op_{\rm GEL} . (t,x(t),x^{(1)} (t) ,\dots ,x^{(n)} (t)) =0 .
\end{equation}

The fractional embedding procedure can be applied to obtain a fractional analogue of this generalized Euler-Lagrange equation. As for the
fractional Euler-Lagrange equation, we need to impose some regularity properties on the Lagrangian.

\begin{defi}
Let $n\in \N^*$ and $\L (t,x,v_1 ,\dots ,v_n )$ be a Lagrangian function of order $n$. The Lagrangian $\L$ is admissible if $\L$ is differentiable
with respect to $x$, holomorphic in each variable $v_i$, $i=1,\dots ,n$, and its restriction to $\R\times \R^d \times (\R ^d )^n$ is real.
\end{defi}

We have:

\begin{thm}
\label{tembgeneulerlagrange}
Let $n\in \N^*$ and $\L$ be an admissible Lagrangian function of order $n$. The $\Emb$-fractional embedding of the generalized Euler-Lagrange
equation associated to $\L$ is given by
$$
\di {\partial L\over \partial x} (t,x(t) ,x^{{}_{\alpha} (1)_{\beta} } (t) ,\dots ,x^{{}_{\alpha} (n)_{\beta} } (t)) +
\di \sum_{i=1}^n (-1)^i \di \left ( \di  \D \right ) ^i
\left [ {\partial L\over \partial v_i} (t,x(t),x^{{}_{\alpha} (1)_{\beta}} (t) ,\dots ,x^{{}_{\alpha} (n)_{\beta} } (t)) \right ] ,
\eqno{(FGEN)^{\alpha ,\beta}_{\mu}}
$$
where
\begin{equation}
x^{{}_{\alpha} (i)_{\beta}} = \left ( \D \right )^i x ,\ \ \ i=1,\dots ,n.
\end{equation}
\end{thm}

The proof is based on the following lemma:

\begin{lem}
\label{emopgel}
Let $n\in \N^*$ and $\L$ be an admissible Lagrangian function of order $n$. The $\Emb$-fractional embedding of the
generalized Euler-Lagrange operator $\Op_{\rm GEL}$ is given by
\begin{equation}
\label{embgeneulerlagrangeop}
\Emb \left ( \Op_{\rm GEL} \right ) =\di\sum_{i=1}^n (-1)^i \left ( \D \right )^i  \circ \di {\partial L\over \partial v_i}
+\di {\partial L\over \partial x} .
\end{equation}
\end{lem}

\begin{proof}
The operator (\ref{geneulerlagrange}) is of the form $\Op_{\mathbf{f}}^{\mathbf{g}}$ with:\\

- $\mathbf{f} =\{ f_i \}_{i=0,\dots ,n}$ where $f_i =(-\mathbf{1} )^i$ and $(-\mathbf{1} )^i (t,x,v_1 ,\dots ,v_n )=(-1)^i$ for $i=0,\dots ,n$.\\

- $\mathbf{g}=\{ g_i \}_{i=0,\dots ,n}$ where $g_0 =\di {\partial L\over \partial x}$ and $g_i =\di {\partial L\over \partial v_i}$, $i=1,\dots ,n$.\\
We then have $\Op_{\rm GEL} =\di \sum_{i=0}^n {}_t F^i \cdot \di {d^i \over dt^i } \circ {}_t G_i$. Using the fractional embedding $\Emb$ of
operators (\ref{formulaembop}) we obtain (\ref{embgeneulerlagrangeop}).
\end{proof}

By definition \ref{defiembequation}, the fractional embedding of equation (\ref{geneulerlagrange}) is given by
\begin{equation}
\Emb \left ( \Op_{\rm GEL} \right ) \cdot (x,\D x ,\dots , \left ( \D \right )^n ) =0 ,
\end{equation}
which reduces to (FGEL) thanks to lemma \ref{emopgel}. This concludes the proof of theorem \ref{tembgeneulerlagrange}.\\

As for the fractional Euler-Lagrange equation, we have an associated {\it coherence problem} (see $\S$.\ref{sectioncoherence}). Next section deals with the
corresponding fractional calculus of variations.

\section{Fractional Ostrogradski formalism}

We develop the notion of extremals and differentiability for functionals depending on higher order, {\it i.e.} power of the fractional operator $\D$.
We prove in particular a coherence theorem.

\subsection{Space of variations and extremals}

The fractional analogue of the space of variations of order $n$ is:

\begin{defi}[Fractional space of variations of order $n$]
Let $n\in \N^*$, we denote by $\varf (n)$ the set of functions $h\in \Ef$ such that
\begin{equation}
(\D )^i h (a)=(\D )^i h (b)=0\ \ \mbox{\rm for}\ \ i=0,\dots ,n-1  .
\end{equation}
\end{defi}

\subsection{The fractional generalized Euler-Lagrange equation}

The explicit form of the differential of a fractional functional of order $n\in \N^*$ is given by the following lemma:

\begin{thm}
\label{difgeneulerlagrange}
Let $n\in \N^*$ and $\L$ be an admissible Lagrangian function of order $n$ with $\Ll$ its associated functional. The functional $\Ll$ is differentiable
at any $x\in \Ef (n)$ and for $h\in \varf (n)$ the differential is given by
\begin{equation}
\label{fracgeneulerlagrange}
d\Ll (x,h)=\int_a^b \left [ \di\sum_{i=1}^n (-1)^i \left ( \Di \right )^i \left [ \di {\partial \L \over \partial v_i} (z_n (t) ) \right ]
+\di {\partial L\over \partial x} (z_n (t)) \right ] dt ,
\end{equation}
where $z_n (t)=\left ( t,x(t),\di\D x(t),\dots ,(\di\D )^n x(t) \right )$.
\end{thm}

\begin{proof}
By linearity of $\D$ we have $\D (x+\epsilon h)=\D x +\epsilon \D h$. We deduce that
\begin{equation}
\Ll (x+\epsilon h)=\di\int_a^b \L (z_n (s) +\epsilon (0,h(s),\D h(s),\dots ,(\D )^n h(s))) \ds .
\end{equation}
Doing a Taylor expansion of $\L (z_n (s) +\epsilon (0,h(s),\D h(s),\dots ,(\D )^n h(s)))$ in $\epsilon$ around $0$, we obtain
\begin{equation}
\label{calculinter1}
\Ll (x+\epsilon h)=\di\int_a^b \left [ \L (z_n (s)) +\epsilon \left (
\di\sum_{i=1}^n \left [ (\D)^i h(s)\right ] {\partial \L\over \partial v_i} (z_n (s))  +\di h(s) {\partial \L \over \partial x} (z_n (s)) \right )
+o(\epsilon ) \right ] \ds .
\end{equation}
Using the product rule (\ref{prodrule}), we obtain by induction the formula
\begin{equation}
\di\int_a^b f(t) (\D )^i h(t) dt =\int_a^b (-1)^i \left [ (\Di )^i f(t) \right ] h(t) dt , \ \ \ i\geq 1.
\end{equation}
Replacing in equation (\ref{calculinter1}) we obtain (\ref{fracgeneulerlagrange}).
\end{proof}

We obtain the fractional analogue of the least-action principle in the fractional Ostrogradski formalism, which follows easily from theorem
\ref{difgeneulerlagrange}.

\begin{thm}[Fractional Ostrogradski least-action principle]
Let $\L$ be an admissible Lagrangian of order $n\in \N^*$, and $\Ll$ its associated functional. A necessary and sufficient condition for a
given function $x\in \Ef (n)$ to be an extremal for $\Ll$ with fixed end points is that it satisfies the fractional generalized Euler-Lagrange equation
$(FGEL)^{\beta ,\alpha}_{-\mu }$.
\end{thm}

We denote by FOLAP the previous derivation of (FGEL).

\subsection{Coherence}

As for the fractional embedding of classical mechanics, the fractional embedding of the Ostrogradski formalism is not always coherent. We have the
usual torsion on the exponents and the parameter $\mu$.

\begin{thm}
Let $\L$ be an admissible Lagrangian of order $n$, $a,b\in \R$, $\alpha ,\beta >0$, then the following diagram commutes
\begin{equation}
\label{diaggeneuler}
\left .
\begin{array}{lll}
\L (t, z_n^{d/dt} (x(t)) ) & \stackrel{\Emb}{\longrightarrow} & \L (t, z_n^{\D} (x(t)) )\\
\mbox{\rm OLAP} \downarrow & & \downarrow \mbox{\rm FOLAP}\\
\mbox{\rm GEL} & \stackrel{\Embi}{\longrightarrow} & \mbox{\rm FGEL}^{\beta ,\alpha }_{-\mu} ,
\end{array}
\right .
\end{equation}
where $z_n^{\mbox{\rm O}} (x(t))=(x(t),\mbox{\rm O} (x(t)) ,\dots ,\mbox{\rm O}^n (x(t))$, with $\mbox{\rm O}$ an operator which can be
$d/dt$ or $\D$.
\end{thm}

Here again, in order to obtain a coherent diagram, we must choose $\mu=0$ and $\alpha=\beta$ at the end, {\it i.e.} on the diagram (\ref{diaggeneuler}).

\section{Nonconservative dynamical systems}

In this section, we provide a unified framework for F. Riewe's approach to nonconservative dynamical systems (\cite{ri1},\cite{ri2}) using the
fractional embedding procedure. For a different approach to this problem we refer to the seminal paper of P.J. Morrison \cite{mori}.

\subsection{Linear friction}

We consider the classical problem of linear friction
\begin{equation}
\label{friction}
-\di {dU \over dx} +\gamma \di {d x\over dt} +m \di {d^2 x \over dt^2} =0 ,
\end{equation}
where $\gamma \in \R$, $m>0$.\\

This system can not be modelled by an ordinary Lagrangian systems. This follows from {\it Bauer theorem} \cite{bau}, which states that it is
impossible to use a variational principle to derive a single linear dissipative equation of motion with constant coefficients.\\

We begin with the following immediate lemma.

\begin{thm}
Let $\L :\R \times \C \times \C \times \C \rightarrow \C$ be the Lagrangian of order $2$ defined by
\begin{equation}
\L (t,x,v_1 ,v_2 )=-U(x)-\di {\gamma \over 2} v_1^2 +\di {1\over 2} m v_2^2 ,
\end{equation}
where $\gamma \in \R$, $m>0$ and $U$ a smooth function.\\

The $\FD^{\alpha}_{\mu}$-fractional embedding of the generalized Euler-Lagrange equation associated to $\L$ is given by
\begin{equation}
-\di {dU\over dx} +\gamma \di (\FD^{\alpha}_{\mu} )^2 x +m (\FD^{\alpha}_{\mu} )^4 x =0 .
\end{equation}
\end{thm}

Specializing to $\mu=-i$, we obtain $\FD^{\alpha}_{-i} =\FD^{\alpha}$. The operator $\FD^{\alpha}$ satisfies a semi-group property (see lemma
\ref{semigroup}). As a consequence, we obtain $\FD^{\alpha} \FD^{\beta} =\FD^{\alpha +\beta}$ for $\alpha >0$, $\beta >0$. We deduce the
following corollary:

\begin{cor}
The $\FD^{\alpha}$-fractional embedding of the generalized Euler-Lagrange equation associated to $\L$ is given by
\begin{equation}
-\di {dU\over dx} +\gamma \di \FD^{2\alpha} x +m \FD^{4\alpha} x =0 .
\end{equation}
\end{cor}

For $\alpha=1/2$ we obtain only fractional operators of integer order, which reduce to ordinary derivatives on the set of smooth functions.
As a consequence, we have proved:

\begin{thm}
The $\FD^{1/2}$-fractional embedding of the generalized Euler-Lagrange equation associated to $\L$ reduce to the equation
\begin{equation}
-\di {dU \over dx} +\gamma \di {d x\over dt} +m \di {d^2 x \over dt^2} =0 ,
\end{equation}
on the set of smooth functions.
\end{thm}

Using the weak-coherence lemma we deduce the stronger result:

\begin{thm}
The solution of the non-conservative systems (\ref{friction}) corresponds to smooth weak-extremals of the $\FD^{1/2}$-fractional functional associated
to $\L$.
\end{thm}

This lemma is not valid for a general Lagrangian systems. However, in our particular case, we have a very simple dependance of the Lagrangian with
respect to velocities. This is this property with the fact that the restriction on $C^2$-functions kills the imaginary part, which implies the
result.

\subsection{The Whittaker equation}

It is in general difficult to see directly if an equation is Lagrangian or not. Following Bateman \cite{bate}, E.T. Whittaker proposes as an answer to a
problem raised by R.C. Tolman (\cite{tol},$\S$.10 $\S$.11) that the equation
\begin{equation}
\label{whit}
\left .
\begin{array}{lll}
\ddot{x} - x & = & 0,\\
\ddot{y} -\dot{x} & = & 0 ,
\end{array}
\right .
\end{equation}
called {\it the Whittaker equation} in what follows, can not be derived by any Lagrangian approach. It seems that up to the author knowledge the question is still {\it open}
today.\\

Following a previous work of F. Riewe \cite{ri2} we give a positive answer to R.C. Tolman's problem in the fractional framework.\\

We begin with the following immediate lemma.

\begin{lem}
Let $\L : \R \times \C^6$ be a Lagrangian of order 2 on $\C^2$ defined by
\begin{equation}
\L (t, x,y,v_1 ,u_1 ,v_2 ,u_2 )=v_2^2 +u_2^2 -u_1 v_1 +x^2 +v_2 y .
\end{equation}
The $\FD_{\mu}^{\alpha}$-embedding of the Euler-Lagrange equation associated to $\L$ is given by
\begin{equation}
\label{fwhitt}
\left .
\begin{array}{lll}
2x(t)+\FD^{\alpha}_{\mu} \left ( \FD^{\alpha}_{\mu} y \right ) -(\FD^{\alpha}_{\mu} )^2  (2 (\FD^{\alpha}_{\mu} )^2 {x}(t)+y(t)) & = & 0 ,\\
\dot{x} (t) +\FD^{\alpha}_{\mu} \left ( \FD^{\alpha}_{\mu} x \right ) - (\FD^{\alpha}_{\mu} )^2 (2 (\FD^{\alpha}_{\mu} )^2 {y}(t)) & = & 0.
\end{array}
\right .
\end{equation}
\end{lem}

Specializing when $\mu=-i$, we obtain $\FD^{\alpha}_{-i} =\FD^{\alpha}$. The operator $\FD^{\alpha}$ satisfies a semi-group property (see lemma
\ref{semigroup}). As a consequence, we have $\FD^{\alpha} \FD^{\alpha} =\FD^{2\alpha}$ for all $\alpha >0$ and equation (\ref{fwhitt}) reduces to:

\begin{thm}
The $\FD^{\alpha}$-embedding of the Euler-Lagrange equation associated to $\L$ is given by
\begin{equation}
\left .
\begin{array}{lll}
2x(t)+\FD^{2\alpha} y  -2\FD^{4\alpha}{x}(t) -\FD^{2\alpha} y(t) & = & 0 ,\\
\dot{x} (t) +\FD^{2\alpha} x - 2\FD^{4\alpha} {y}(t) & = & 0.
\end{array}
\right .
\end{equation}
\end{thm}

Choosing $\alpha=1/2$, we obtain fractional operators of integer order. On the set of sufficiently smooth functions, these operators reduces to the
ordinary derivatives. As a consequence, we have proved:

\begin{thm}
The restriction of the $\FD^{1/2}$-fractional embedding of the Euler-Lagrange equation associated to $\L$ is given by the Whittaker equation.
\end{thm}

Moreover, using the weak least action principle, we have the stronger result:

\begin{thm}
Solutions of the Whittaker equation correspond to smooth weak-extremals of the $\FD^{1/2}$-fractional functional associated to $\L$.
\end{thm}

\newpage
\part{Frational embedding of continuous Lagrangian systems}
\label{part8}
\setcounter{section}{0}

The wave equation is an example of a partial differential equation which can be derived via a variational principle. In that case, the functional
depends on several variables (at least two). Such functional arise in mechanical problems involving systems with infinitely many degrees of freedom
like strings, membrane, etc. These systems are known as {\it continuous Lagrangian systems}. We use the fractional embedding procedure to define the
fractional analogue of these equations. In particular, we obtain a fractional Lagrangian formulation of the fractional wave equation introduced by
Schneider and Wyss \cite{sw}. This suggest also a canonical way to generalize to the fractional framework a given differential equation: if the
underlying equation possesses an additional structure like the Lagrangian case, one must provide an extension which respect this structure. The
fractional embedding procedure is well adapted to deal with such an extension. We also derive a fractional Lagrangian formulation for the fractional
diffusion equation.

\section{Continuous Lagrangian systems}

Let $d\in \N$. We consider a Lagrangian function $\L$ defined on $\R \times \R^d \times \R \times \C \times \R^d$ and denoted by
\begin{equation}
\label{density}
\L (t,x_1 ,\dots ,x_d ,y,v,w_1 ,\dots ,w_d ) .
\end{equation}
In what follows, we use the terminology of {\it Lagrangian density} for a function $\L$ of the form (\ref{density}). We denote
$x$ for $(x_1 ,\dots ,x_d)$.\\

Let $\Re$ be a fixed region of $\R^d$ and $a <b$, $a,b\in \R$. We consider the functional
\begin{equation}
\label{extfield}
{\cal L}_{a,b,\Re} (u) =\di \int_a^b \int_{\Re} \L (t,x,u(t,x),\partial_t u (t,x) ,\partial_x u (t,x) ) \ dx \, dt ,
\end{equation}
acting on a function $u:\R\times \R^d \longrightarrow \R$ which is usually called a {\it field}, which is of class $C^1$ in all its variables and where
\begin{equation}
\partial_x u(t,x)=(\partial_{x_1} u(t,x) ,\dots ,\partial_{x_d} u(t,x)) .
\end{equation}
A {\it variation} for a field $u(t,x)$ is defined as a function of the form
\begin{equation}
u_{\epsilon} (t,x)=u(t,x) +\epsilon h(t,x) ,
\end{equation}
with $h$ satisfying the {\it boundaries conditions}
\begin{equation}
h(a,x)=h(b,x)=0\ \ \ \mbox{\rm and}\ \ \ h(t,\partial \Re )=0 ,
\end{equation}
where $\partial \Re$ denotes the boundary of $\Re$.

\begin{thm}[Euler-Lagrange equation for fields]
Extremals of the functional (\ref{extfield}) are solutions of the Euler-Lagrange equation for fields (ELF):
\begin{equation}
\di {\partial \L\over \partial y} -\di {\partial \over \partial t} \left ( \di {\partial \L \over \partial v} \right ) -\di\sum_{i=1}^d
\di {\partial \over \partial x_i} \left [ \di {\partial \L \over \partial w_i} \right ] =0 .
\end{equation}
\end{thm}

We refer to (\cite{gf} Chapter 7) for more details.\\

As in the previous parts, we associate to (ELF) the differential operator
\begin{equation}
\Op_{\mbox{\rm (ELF)}} =\di {\partial \L\over \partial y} -\di {\partial \over \partial t} \circ {\partial \L \over \partial v} -\di\sum_{i=1}^d
\di {\partial \over \partial x_i} \circ \di {\partial \L\over \partial w_i} ,
\end{equation}
acting on $2d+2$-uplet $(t,x,u(t,x),\partial_t u (t,x), \partial_{x_1} u (t,x) ,\dots ,\partial_{x_d} u (t,x) )$.\\

Using the fractional embedding procedure we obtain a fractional Euler-Lagrange equation for fields.

\begin{thm}
Let $\L$ be an admissible Lagrangian density. The $\Emb$-fractional embedding of the Euler-Lagrange equation for fields (ELF) associated to
$\L$ is given by
$$
\di {\partial \L\over \partial y} (z_{\alpha ,\beta} (t,x)) -\di \D \left [ {\partial \L \over \partial v} (z_{\alpha ,\beta} (t,x) \right ]-
\di\sum_{i=1}^d \di {\partial \over \partial x_i} \left [ {\partial \L\over \partial w_i} (z_{\alpha ,\beta} (t,x)) \right ]
\eqno{(FELF)^{\alpha ,\beta}_{\mu}}
$$
where $z_{\alpha ,\beta} (t,x) =(t,x,u(t,x),\D u(t,x) ,\partial_{x_1} u (t,x) ,\dots ,\partial_{x_d} u(t,x))$.
\end{thm}

\begin{proof}
The differential operator $\Op_{\rm (ELF)}$ is considered as a time-differential operator, {\it i.e.} that we consider the operator
$\Op_{\mathbf{f}}^{\mathbf{g}}$ with
\begin{equation}
\mathbf{f}=\left ( \di {\partial \L \over \partial y} -\di\sum_{i=1}^d \di {\partial\over \partial x_i} \circ \di {\partial \L \over \partial w_i} ,\mathbf{1}
\right ) ,
\end{equation}
and
\begin{equation}
\mathbf{g} =\left ( \mathbf{1} ,-\di {\partial \L \over \partial v} \right ) .
\end{equation}
Using definition \ref{defiembequation}, we obtain (\mbox{\rm FEL}). This concludes the proof.
\end{proof}

The previous theorem can be generalized in various ways. In particular, one can define a generalized embedding procedure, assuming that not only time
derivatives, but also spatial derivatives $\partial/\partial x_i$, $i=1,\dots ,d$, are replaced by a fractional derivatives. In this case however, we
have some technical difficulties (see remark \ref{remarkgreen}). This will be discussed in another work.\\

As for the previous cases, we develop the associated fractional calculus of variations.

\section{Fractional continuous lagrangian systems}

Using the fractional embedding procedure, we look for the following class of {\it fractional densities}:

\begin{defi}
Let $\L$ be an admissible Lagrangien density. The fractional functional associated to $\L$ is defined by
\begin{equation}
\label{functionalcontinuous}
\Llf (u) =\di \int_a^b \int_{\Re} \L (t,x,u(t,x),\di\D u (t,x) ,\partial_x u (t,x) ) \ dx \, dt ,
\end{equation}
for fields $u(t,x) \in {}^{\alpha}_a \F^{\beta}_b (\Re )$, the set of fields smooth with respect to $x$ and in $\Ef$ with respect to $t$.
\end{defi}

We consider two spaces of variations for fields:

\begin{defi}[Spaces of variations for fields]
We denote by $\Var^{\alpha ,\beta ,\Re} (a,b)$ the set of fields satisfying
\begin{equation}
\Var^{\alpha ,\beta} (a,b ,\Re )=\{ h(t,x) , h_t \in C^1 ,\ h_x \in \Ef ,\ h(a,x)=h(b,x)=0,\ h(t,\partial \Re )=0 \} ,
\end{equation}
and by $\Var^{\alpha}_0 (a,b)$ the set of fields defined by
\begin{equation}
\Var^{\alpha}_0 (a,b ,\Re )=
\left \{
\begin{array}{c}
h(t,x) , h_t \in C^1 ,\ h_x \in \Ef ,\ h(a,x)=h(b,x)=0,\\
\ h(t,\partial \Re )=0 ,\ \adta h=\tdba h
\end{array}
\right \} .
\end{equation}
\end{defi}

As usual the condition $\adta h=\tdba h$ allows us to obtain a symmetric product rule which is fundamental for the coherence problem.

\begin{defi}
Let $u(t,x)$ be a given field. A $\P$-variation of $u$ is a one-parameter $\epsilon \in \R$ family of fields $u_{\epsilon} (t,x)$ defined by
\begin{equation}
u_{\epsilon} (t,x)=u(t,x)+\epsilon h(t,x) ,\ \ \ h\in \P,
\end{equation}
where $\P$ can be $\Var^{\alpha ,\beta} (a,b)$ or $\Var^{\alpha}_0 (a,b)$.
\end{defi}

We denote by $\Ll_{a,b}$ the functional $\Ll_{a,b}^{\alpha ,\beta}$ when $\P =\Var^{\alpha ,\beta} (a,b ,\Re )$ and $\Ll_{a,b}^{\alpha ,\alpha}$
when $\P =\Var^{\alpha} (a,b,\Re )$.

\begin{defi}
Let $\L$ be an admissible Lagrangian density and $\Ll_{a,b}$ the associated functional. The functional $\Ll_{a,b}$ is called $\P$-differentiable
at $u$, where $u$ is a field, if
\begin{equation}
\Ll_{a,b} (u+\epsilon h ) -\Ll_{a,b} (u)=\epsilon d\Ll_{a,b} (u,h) +o(\epsilon ) ,
\end{equation}
for all $h\in \P$, where $d\Ll_{a,b} (u,h)$ is a linear functional of $h$.
\end{defi}

The differential of the fractional functional (\ref{functionalcontinuous}) is given by:

\begin{thm}
Let $\L$ be an admissible Lagrangian density and $\Llf$ the associated fractional functional. The functional $\Llf$ is differentiable
at any fields $u\F$ and for all $h\in \Var^{\alpha ,\beta} (a,b,\Re )$ the differential is given by
\begin{equation}
\di {\partial \L\over \partial y} (z_{\alpha ,\beta} (t,x)) -\di \Di \left [ {\partial \L \over \partial v} (z_{\alpha ,\beta} (t,x) \right ] +
\di\sum_{i=1}^d \di {\partial \over \partial x_i} \left [ {\partial \L\over \partial w_i} (z_{\alpha ,\beta} (t,x)) \right ] ,
\end{equation}
where $z_{\alpha ,\beta} (t,x) =(t,x,u(t,x),\D u(t,x) ,\partial_{x_1} u (t,x) ,\dots ,\partial_{x_d} u(t,x))$.
\end{thm}

\begin{proof}
As the left and right (RL) derivatives are linear operators, we have
\begin{equation}
\D (u+\epsilon h )= \D u +\epsilon \D h .
\end{equation}
As a consequence, we obtain
\begin{equation}
\Llf (u+\epsilon h)=\di\int_a^b \int_{\Re} \L ( z_u (t,x)+\epsilon z_h (t,x)) \ \ dx\, dt ,
\end{equation}
where for a field $f(t,x)$ we denote
\begin{equation}
z_f (t,x)=(t,x,f(t,x),\D f (t,x) ,\partial_{x_1} f(t,x) ,\dots ,\partial_{x_d} f (t,x)) .
\end{equation}
Doing a Taylor expansion of $\L ( z_u (t,x)+\epsilon z_h (t,x))$ in $\epsilon$ around $0$, we obtain
\begin{equation}
\left .
\begin{array}{lll}
\L ( z_u (t,x)+\epsilon z_h (t,x)) & = & \L (z_u (t,x)) \\
 & & +\epsilon \left [
\di {\partial \L \over \partial y} (z_u (t,x)) \dot z_h (t) + \D h \di {\partial \L \over \partial v} (z_u (t,x)) +
\di\sum_{i=1}^d \di {\partial \L \over \partial w_i} (z_u (t,x)) \partial_{x_i} h \right ]\\
& &  +o(\epsilon ).
\end{array}
\right .
\end{equation}
Using the product rule (\ref{prodrule}) we obtain
\begin{equation}
\di \int_a^b \int_{\Re } \D h \di {\partial \L \over \partial v} (z_u (t,x)) dx dt =-
\di \int_a^b \int_{\Re } h \Di \left [ \di {\partial \L \over \partial v} (z_u (t,x)) \right ] dx dt
\end{equation}
as long as $h(a,x)=h(b,x)=0$.\\

Moreover, using the multidimensional Green theorem (see \cite{wi} p.223), we have
\begin{equation}
\label{form19}
\int_a^b \int_{\Re} \di\sum_{i=1}^d \di {\partial \L \over \partial w_i} (z_u (t,x)) \partial_{x_i} h \ dx\, dt=
\int_a^b \int_{\Re} \sum_{i=1}^d \di {\partial \over \partial x_i} \left [ \di {\partial \L \over \partial w_i} (z_u (t,x)) \right ] h\ dx\,  dt,
\end{equation}
as long as $h(t,\partial \Re )=0$. This concludes the proof.
\end{proof}

\begin{rema}
\label{remarkgreen}
In the global fractional framework, {\it i.e.} when we use fractional derivatives also in the $x$ variables, we see from formula (\ref{form19}) that we will need an analogue
of the multidimensional Green theorem in the fractional case.
\end{rema}

We denote by ${\cal D}^{\alpha}_{\mu}$ the operator ${\cal D}^{\alpha ,\alpha }_{\mu}$ and $\F^{\alpha} (\Re )$ the set
${}^{\alpha}_a \F^{\alpha}_b (\Re )$.

\begin{thm}[Real variations]
Let $\L$ be an admissible Lagrangian density and $\Llf$ the associated fractional functional. The functional $\Llf$ is differentiable
at any fields $u\in \F^{\alpha} (\Re ) $ and for all $h\in \Var^{\alpha} (a,b,\Re )$ the differential is given by
\begin{equation}
\di {\partial \L\over \partial y} (z_{\alpha } (t,x)) -\di {\cal D}^{\alpha}_{\mu} \left [ {\partial \L \over \partial v} (z_{\alpha } (t,x) \right ]-
\di\sum_{i=1}^d \di {\partial \over \partial x_i} \left [ {\partial \L\over \partial w_i} (z_{\alpha } (t,x)) \right ] ,
\end{equation}
where $z_{\alpha } (t,x) =(t,x,u(t,x),{\cal D}^{\alpha}_{\mu} u(t,x) ,\partial_{x_1} u (t,x) ,\dots ,\partial_{x_d} u(t,x))$.
\end{thm}

The proof follows the same lines.\\

An immediate consequence of this lemma is:

\begin{thm}[Fractional least-action principle for fields]
Let $\L$ be an admissible Lagrangian density and $\Llf$ the associated functional. A necessary and sufficient condition for a field $u$
to be a $\Var^{\alpha ,\beta} (a,b ,\Re )$-extremal is that it satisfies the fractional Euler-Lagrange equation for fields $(FELF)^{\beta ,\alpha}_{-\mu}$.
\end{thm}

We denote by ${\cal L}^{\alpha}_{a,b}$ the functional ${\cal L}^{\alpha ,\alpha}_{a,b}$. Using the set of real variations, we have:

\begin{thm}[Weak Fractional least-action principle for fields]
Let $\L$ be an admissible Lagrangian density and ${\cal L}^{\alpha}_{a,b}$ the associated functional. A necessary and sufficient condition for a field $u$
to be a $\Var^{\alpha} (a,b\Re )$-extremal is that it satisfies the fractional Euler-Lagrange equation for fields $(FELF)^{\alpha}_{\mu}$.
\end{thm}

\section{The fractional wave equation}

In this section, we derive the {\it fractional wave equation} defined by Schneider and Wyss \cite{sw} as the extremals of a fractional continuous
Lagrangian systems.\\

The equation describing waves propagating on a stretched string of constant linear mass density $\rho$ under constant tension $T$ is
\begin{equation}
\rho \di {\partial^2 u (t,x)\over \partial t^2} =T \di {\partial^2 u(t,x)\over \partial x^2} ,
\end{equation}
where $u(t,x)$ denotes the amplitude of the wave at position $x$ along the string at time $t$. The wave equation corresponds to the extremals of the
generalized functional associated to the Lagrangian systems
\begin{equation}
\label{wave}
\L (t,x,y,v,w)= \di {\rho \over 2} v^2 - \di {T\over 2} w^2 .
\end{equation}

In \cite{sw}, the authors define the fractional analogue of the wave equation by changing the classical derivative by a fractional one. Using our
notations, the definition of the fractional wave equation is:

\begin{defi}
The fractional wave equation of order $\alpha$ is the fractional differential equation
\begin{equation}
-\rho \FD^{2\alpha} u =T\di {\partial^2 u\over\partial x^2} .
\end{equation}
\end{defi}

A natural demand with respect to this generalization which is just a formal manipulation on equations, is to keep a more structural property of the
wave equation, namely the fact that it derives from a least-action principle. Using our fractional embedding procedure, we are able to explicit such a
fractional Lagrangian framework for the fractional wave equation.\\

In the following we work with the fractional embedding associated to $\mathbf{D}^{\alpha}_{\mu}$ which corresponds to $\D$ when
$\beta=\alpha$ and $a=-\infty$, $b=+\infty$.

\begin{thm}
The $\FD^{\alpha}_{\mu}$-fractional embedding of the continuous Euler-Lagrange equation associated to (\ref{wave}) is given by
\begin{equation}
-\rho \mathbf{D}^{\alpha}_{\mu} \circ \mathbf{D}^{\alpha}_{\mu}  u =T \partial^2_{x^2} u .
\end{equation}
\end{thm}

We can specialized by choosing $\mu=-i$. In that case $\mathbf{D}^{\alpha}_{-i} =\FD^{\alpha}$ and satisfies a semi-group property (see
lemma \ref{semigroup}). As a consequence, we obtain:

\begin{cor}
The $\FD^{\alpha}$-fractional embedding of the continuous Euler-Lagrange equation associated to (\ref{wave}) is given by
\begin{equation}
\label{fwave}
-\rho \FD^{2\alpha}  u =T \partial^2_{x^2} u .
\end{equation}
\end{cor}

Moreover, using the weak coherence theorem, we have:

\begin{thm}
Solutions of the fractional wave equation (\ref{fwave}) of order $\alpha$ corresponds to weak-extremals of the $\FD^{\alpha}$-fractional functional
associated to $\L$.
\end{thm}

Up to the author knowledge, this is the first time that the fractional wave equation is derived via a fractional variational principle.
In particular, the previous derivation has the advantage to keep the continuous Lagrangian structure underlying the classical wave equation.

\section{A remark on fractional/classical equations}

The previous remark can be uses as a conceptual guideline to generalize classical equations of physics in the fractional framework. If the
classical equation possesses an additional structure, for example Lagrangian, then we must extend this equation keeping this additional
structure, generalized in a natural way. The main remark is that equations by themselves do not have a universal significance their form depending mostly
on the coordinates systems being used to derive them. On the contrary the underlying first principle like the least-action principle carry an
information which is of physical interest and not related to the coordinates system which is used. At least this point of view explain the importance of
coherence theorems in all the existing embedding theories of dynamical systems.

\section{The fractional diffusion equation}

The {\it fractional diffusion equation} of order $0<\alpha <1$ is defined by
\begin{equation}
\label{fracdiffusion}
\FD^{\alpha} u(t,x) =a^2 \di {\partial^2 u(t,x)\over \partial x^2} .
\end{equation}
It is defined by Wyss \cite{wy}. For $\alpha =1$ we recover the classical diffusion equation.\\

The aim of this section is to derive a fractional Lagrangian formulation of the fractional wave equation of order $0<\alpha \leq 1$, then including the
classical diffusion equation. In the contrary of the fractional wave equation, the diffusion equation is recovered thanks to a still fractional
variational principle.\\

Let us consider the Lagrangian function $\L$ defined on $\R\times \R^2 \times \C \times \C$ by
\begin{equation}
\label{diffusionL}
\L (t,x,y,v,w)=\di {1\over 2} v^2 - \di {a^2 \over 2} w^2 ,
\end{equation}
where $\rho \in \R$, $a\in \R$.\\

As for the fractional wave equation, we denote by $\mathbf{D}_{\mu}^{\alpha}$ the quantity $\D$ when $\beta=\alpha$, $a=-\infty$ and
$b=+\infty$, {\it i.e.} working with the fractional derivatives $\FD$ anf $\FD_*$.

\begin{thm}
The $\mathbf{D}^{\alpha /2}_{\mu}$-fractional embedding of the continuous Euler-Lagrange equation associated to (\ref{diffusionL}) is given by
\begin{equation}
\mathbf{D}^{\alpha /2}_{\mu} \circ \mathbf{D}^{\alpha/2}_{\mu} u =a^2 \di {\partial^2 u \over \partial x^2} .
\end{equation}
\end{thm}

Choosing $\mu=-i$, we obtain $\mathbf{D}^{\alpha}_{-i} =\FD^{\alpha}$. As $\FD^{\alpha}$ satisfies a semi-group property, we finally obtain:

\begin{thm}
The $\FD^{\alpha/2}$-fractional embedding of the continuous Euler-Lagrange equation associated to (\ref{diffusionL}) is given by
\begin{equation}
\FD^{\alpha} u =a^2 \di {\partial^2 u \over \partial x^2} .
\end{equation}
\end{thm}

It must be noted that even for $\alpha=1$, the diffusion equation is recovered using a fractional embedding procedure, namely the $\FD^{1/2}$-fractional
embedding procedure.\\

The main result of this section is that this fractional embedding of the diffusion equation has an additional structure, a Lagrangian one.

\begin{thm}
Solutions of the fractional wave equation (\ref{fracdiffusion}) of order $0<\alpha \leq 1$ corresponds to weak-extremals of the $\FD^{\alpha /2}$-fractional functional
associated to (\ref{diffusionL}).
\end{thm}

This result seems new, even for the case $\alpha=1$.

\newpage
\part*{Conclusion and perspectives}

Part of this paper can be generalized in various ways.\\

i) We must extend the fractional embedding procedure in order to cover partial differential equations. In particular, the general fractional
reaction-diffusion equation must be studied in this setting.\\

ii) Of special interest is the Schr\"odinger equation. With F. Ben Adda (\cite{bc2},\cite{bc3}) we have obtain this equation following an idea of L. Nottale
\cite{no} and results about the local fractional calculus \cite{bc1}. A problem is to try to do the same computations using left and
right RL fractional derivatives. The main difficulty is the absence of a chain rule for these operators. In that case, one must obtain a one-parameter
family of partial differential equations, depending on the non-integer order $\alpha$ of differentiation. Following \cite{cr3} the idea is to
prove that the Schr\"odinger equation can be recovered only for $\alpha =1/2$ which gives strong constraints on the underlying nature of space-time in
Nottale's framework of the {\it Scale relativity}.\\

iii) The previous idea can in fact be generalized to stochastic processes following our previous work with S. Darses (\cite{cd1} \cite{cd3}). In that case,
one want to consider fractional Brownian processes via the stochastic embedding procedure. The stochastic derivative defined in \cite{cd3} must be generalized and
the paper by S. Darses and I. Nourdin \cite{dn} can be considered as a first step in this direction.\\

iv) We have use the left and right RL fractional derivatives, but the same scheme can be used with another fractional calculus leading to a different
fractional embedding theory. This depends mostly on the type of applications one want to consider. Of particular interest seems to be the Caputo fractional
derivatives \cite{skm}. A forthcoming paper will be devoted to this case.

\newpage
${}$
\newpage
\part*{Notations}
\vskip 1cm
\noindent {\bf Operators}
\begin{itemize}
\item $\adta$: left Riemann-Liouville fractional derivative of order $\alpha$
\item $\tdba$: right Riemann-Liouville fractional derivative of order $\alpha$
\item $\FD^{\alpha}$: left fractional derivative of order $\alpha$
\item $\FD_*^{\alpha}$: right fractional derivative of order $\alpha$
\item $\D$: fractional operator of order $(\alpha ,\beta )$ and parameter $\mu$
\item $\FD^{\alpha ,\beta}_{\mu}$: fractional operator of order $(\alpha ,\beta )$ and parameter $\mu$ associated to $\FD$ and $\FD_*$
\item $\mbox{\rm d}^{\alpha}$: the reversible fractional derivative corresponding to $\D$ for $\beta=\alpha$, $\mu=0$.
\end{itemize}
\vskip 6mm
{\bf Functional spaces}
\begin{itemize}
\item ${}_a^{\alpha} \mbox{\rm E}$: left fractional Riemann-Liouville derivative space
\item $\mbox{\rm E}_b^{\alpha}$: right fractional Riemann-Liouville derivative space
\item ${}_a^{\alpha} \mbox{\rm E}_b^{\beta}$: functional space associated to the fractional operator $\D$
\item $AC([a,b])$: absolutely continuous functions on $[a,b]$
\item $C_0^{\infty} (I)$: functions of $C^{\infty} (I)$ vanishing ouside a compact subset $K$ of $I$.
\item $J_L^{\alpha} (\R )$: left fractional derivative space
\item $J_R^{\alpha} (\R )$: right fractional derivative space
\item $J_{L,0}^{\alpha} (I)$, $J_{R,0}^{\alpha} (I)$: closure of $C_0^{\infty} (I)$ under their respective norms
\item $H_0^{\alpha} (I)$: fractional Sobolev space of order $\alpha$
\end{itemize}
\vskip 6mm
{\bf Operators and objects associated to the fractional embedding procedure}
\begin{itemize}
\item $\Emb$: the fractional embedding procedure associated to $\D$
\item $\Rev$: the reversible fractional embedding procedure, {\it i.e.} ${}_a^{\alpha} \mbox{\rm Emb}^{\alpha}_b (0)$
\item $\Op_{\rm (EL)}$: the Euler-Lagrange operator
\item $\rm EL$: the Euler-Lagrange equation
\item $\mbox{\rm FEL}^{\mu}_{\alpha ,\beta}$: the fractional Euler-Lagrange equation
\item $\mbox{\rm\bf Var}^{\alpha}$: the real space of variations
\item ${\cal L}^{\alpha ,\beta}_{a,b}$: the fractional functional associated to $\L$
\end{itemize}
\newpage
${}$


\newpage
\begin{thebibliography}{15}
\bibitem{agra} Agrawal O.P., Formulation of Euler-Lagrange equations for fractional variational problems, J. Math. Anal. Appl. 272 (2002) 368-379.
\bibitem{and} Andler M., Jean Leray (1906-1998), Proceeding of the American Phylosophical Society Vol. 144, no. 4 (2000).
\bibitem{ar} Arnold V.I., {\it Mathematical methods of classical mechanics}, 2d Edition, Springer, 1989.
\bibitem{av} Avez A., {\it Calcul diff\'erentiel}, Masson, 1983.
\bibitem{bate} Bateman H., On dissipative systems and related variational problems, Phys. Rev. {\bf 38}, 815-819 (1931).
\bibitem{bc1} Ben Adda F., Cresson J., About non-differentiable functions, J. Math. Anal. Appl. 263 (2001) 721-737.
\bibitem{bc2} Ben Adda F., Cresson J., Fractional differential equations and the Schr\"odinger equation, Appl. Math. Comp. 161 (2005) 323-345.
\bibitem{bc3} Ben Adda F., Cresson J., Quantum derivatives and the Schr\"odinger equation, Chaos Solitons and Fractals 19 (2004)
1323-1334.
\bibitem{bau} Bauer P.S., Proc. Nat. Acad. Sci. USA {\bf 17}, 311 (1931).
\bibitem{clz} Carreras B.A., Lynch V.E., Zaslavsky G.M., Anomalous diffusion and exit time distribution of particle tracers in plasma turbulence models,
Phys. Plasmas 8(12) 5096-5103 (2001).
\bibitem{cm} Carpinteri A., Mainardi F., {\it Fractals and fractional calculus in continuum mechanics}, Springer-Verlag, Wien, 1997.
\bibitem{comte} Comte F., Op\'erateurs fractionnaires en Econom\'etrie et en Finance, Pr\'epublication MAP5 2001-3.
\bibitem{com} Compte A., Stochastic foundations of fractional dynamics, Phys. Rev. E {\bf 53} 4 (1996) 4191-4193.
\bibitem{cr1} Cresson J., Non-differentiable variational principles, J. Math. Anal. Appl. 307 (2005) 48-64.
\bibitem{cr2} Cresson J., {\it Th\'eories de plongement des syst\`emes dynamiques - Un programme}, 21.p, 2005.
\bibitem{cr3} Cresson J., Scale calculus and the Schr\"odinger equation, J. Math. Phys. 44 (2003) 4907-4938.
\bibitem{cr4} Cresson J., Quantum embedding and partial differential equations, {\it in preparation}, 2006.
\bibitem{cd1} Cresson J., Darses S., Plongement stochastique des syst\`emes dynamiques, C. R. Acad. Sci. Paris. Ser. I 342 (2006) 333-336.
\bibitem{cd2} Cresson J., Darses S., Th\'eor\`eme de Noether stochastique, Pr\'epublication de l'I.H.\'E.S. 06/25 (2006)  7.p.
\bibitem{cd3} Cresson J., Darses S., {\it Stochastic embedding of dynamical systems}, Pr\'epublications de l'I.H.\'E.S. 06/27 (2006) 92.p.
\bibitem{cd4} Cresson J., Darses S., Stochastic perturbation theory, {\it in preparation}, 2006.
\bibitem{dn} Darses S., Nourdin I., Stochastic derivatives for fractional diffusions, Pr\'epublications Paris VI (2006).
\bibitem{er} Erwin V.J., Roop J.P., Variational formulation for the stationary fractional advection dispersion equation,
Numer. Meth. P.D.E., 22 pp. 558-576, (2006).
\bibitem{gasto} Frederico G., Torres D., A formulation of Noether's theorem for fractional problems of the calculus of variations, preprint (2005).
\bibitem{gf} Gelfand I.M., Fomin S.V., {\it Calculus of variations}, Dover Publication Inc, 2000.
\bibitem{gm} Gorenflo R., Mainardi F., Fractional calculus: integral and differential equations of fractional order, in \cite{cm} 223-276.
\bibitem{hw} Henry B.I., Wearne S.L., Fractional reaction-diffusion, Physica A 276 (2000) 448-455.
\bibitem{sch} Kosmann-Scwarzback Y., {\it Les th\'eor\`emes de Noether; invariance et lois de conservation au XX \`eme si\'ecle}, Les \'editions
de l'\'Ecole Polytechnique (2004).
\bibitem{ler} Leray J., Sur le mouvement d'un liquide visqueux emplissant l'espace. Acta Math. 63, 193-248 (1934).
\bibitem{marmi} Marmi S., Chaotic behavior in the Solar system [following J. Laskar], S\'eminaire Bourbaki 51 \`eme ann\'ee, 1998-1999? no. 854.
\bibitem{mnn} Le Mehaute A., R.R. Nigmatullin, R. Nivanen, {\it Fl\`eches du temps et g\'eom\'etries fractales}, Ed. Hermes, Coll. Syst\`emes Complexes,
1998.
\bibitem{mr} Miller K.S., Ross B., {\it An introduction to fractional calculus and fractional differential equations}, John Wiley and Sons, New York,
1993.
\bibitem{mori} Morrison P.J., A paradigm for joined Hamiltonian and dissipative systems, Physica D 18 (1986) 410-419.
\bibitem{muba} Muslih S.I., Baleanu D., Formulation of Hamiltonian equations for fractional variational problems, arXiv:math-ph/0510029 (2006).
\bibitem{no} Nottale L., {\it Fractal space-time and microphysics}, World Scientific, 1993.
\bibitem{os} Oldham K.B., Spanier J., {\it The fractional calculus}, Academic Press, New York, 1974.
\bibitem{pod} Podlubny I., {\it Fractional differential equations}, Academic Press, New York, 1999.
\bibitem{ri1} Riewe F., Nonconservative Lagrangian and Hamiltonian mechanics, Phys. Rev. E {\bf 53} 2 (1996) 1890-1899.
\bibitem{ri2} Riewe F., Mechanics with fractional derivatives, Phys. Rev. E {\bf 55} 3 (1997) 3581-3592.
\bibitem{skm} Samko S.G., Kilbas A.A., Marichev O.I., {\it Fractional integrals and derivatives: Theory and applications}, Gordon and Breach, New York,
1998.
\bibitem{sw} Schneider W.R., Wyss W., Fractional diffusion and wave equations, J. Math. Phys. 30 (1989) 134-144.
\bibitem{sww} Shlesinger M.F., West B.J., Klafter J., L\'evy dynamics of enhanced diffusion: application to turbulence, Phys. Rev. Lett 58 (11) 1100-1103
(1987).
\bibitem{tol} Tolman R.C., {\it Statistical mechanics}, Chemical Catalog Co. NewYork 1927.
\bibitem{wi} Widder D.V., {\it Advanced calculus} 2d edition, Dover, N.Y., 1961.
\bibitem{wy} Wyss W., The fractional diffusion equation, J. Math. Phys. 27 (1986) 2782-2785.
\bibitem{zu} Zaslavsky G.M., {\it Hamiltonian chaos and fractional dynamics}, Oxford University Press (2005).
\end{thebibliography}
\end{document}